\documentclass[12pt,reqno]{amsart}

\usepackage{fullpage,amsmath,amssymb,amsfonts,mathrsfs,bm,graphicx,hyperref}

\def\bbone{{\mathchoice {\rm 1\mskip-4mu l} {\rm 1\mskip-4mu l}
{\rm 1\mskip-4.5mu l} {\rm 1\mskip-5mu l}}}

\newtheorem{theorem}{Theorem}[section]
\newtheorem{remark}{Remark}[section]

\newtheorem{corollary}{Corollary}[section]
\newtheorem{lemma}{Lemma}[section]
\newtheorem{proposition}{Proposition}[section]

\begin{document}

\author{Abdelmalek Abdesselam}
\address{Abdelmalek Abdesselam, Department of Mathematics,
P. O. Box 400137,
University of Virginia,
Charlottesville, VA 22904-4137, USA}
\email{malek@virginia.edu}

\title{Proof of a conjecture by Starr and log-concavity for random commuting permutations}

\begin{abstract}
We prove a conjecture by Shannon Starr regarding the asymptotics for the number of tuples of commuting permutations with given number of joint orbits. These numbers generalize unsigned Stirling numbers of the first kind which count how many single permutations have a given number of cycles. In the case of pairs of permutations, these numbers are related to D'Arcais polynomials and the Nekrasov-Okounkov formula.
As a consequence of the above asymptotics, we confirm a log-concavity conjecture in the regime of typical values for the number of joint orbits.
As a result of possible indepentent interest in applied mathematics and mathematical physics, we also provide detailed asymptotics, using Mellin transform techniques, for certain multiple series or multivariate Ramanujan sums which are related to ordinary generating functions of Dirichlet convolutions of power laws. Besides these multiple sums asymptotics, our proofs use bivariate saddle point analysis related to the Meinardus theorem in the delicate case of multiple poles for the associated Dirichlet series. 
\end{abstract}

\maketitle


\section{Introduction}

\subsection{Main results}

For any $\ell\in\mathbb{Z}_{>0}:=\{1,2,\ldots\}$, any $n\in\mathbb{N}:=\{0,1,2,\ldots\}$, and any integer $k$ with $0\le k\le n$, we define the enumerative quantity $A(\ell,n,k)$ as follows. For an ordered $\ell$-tuple $(\sigma_1,\ldots,\sigma_{\ell})$
of permutations in the symmetric group $\mathfrak{S}_n$, we define a joint orbit of such a tuple as an orbit in the set $[n]:=\{1,2,\ldots,n\}$ for the natural permutation action of the subgroup of $\mathfrak{S}_n$ generated by the permutations $\sigma_1,\ldots,\sigma_{\ell}$. By definition, $A(\ell,n,k)$ is the number of tuples $(\sigma_1,\ldots,\sigma_{\ell})$ which have exactly $k$ joint orbits. Note that when $\ell=1$, i.e., the case of a single permutation, these orbits simply are the cycles of that permutation. For higher $\ell$, these look like discrete $\ell$-dimensional tori inside the set $[n]$ (see~\cite{AbdesselamBDV} and the appendix of~\cite{Starr}).

For $\ell,n,k$ as above, we define the quantity $\mathscr{S}_{\ell}(n,k)$, introduced by Shannon Starr in~\cite[Conj. 2.3]{Starr}, by
\[
\mathscr{S}_{\ell}(n,k):=
-\ell\ k\ln k+k\left(\ 
(\ell-1)\ln n-\ell\ \ln(\ell-1)+\ell+\ln \mathcal{K}_{\ell}
\ \right)
\]
where, in keeping with the notation of~\cite{AbdesselamS},
\[
\mathcal{K}_{\ell}:=(\ell-1)!\ \zeta(2)\zeta(3)\cdots\zeta(\ell)\ ,
\]
in terms of special values of the Riemann zeta function.
Our first main result is the following theorem.

\begin{theorem}\label{mainthm1}
For any $\ell\ge 4$, there exists a sequence of explicit constants $H_1,H_2,\ldots,H_{\ell-1}$, such that for all $s\in(0,1)$, and all sequences $(k_n)_{n\ge 1}$, with
$1\le k_n\le n$ that are such that
$k_n\sim s n^{\frac{\ell-1}{\ell}}$, we have the asymptotic equivalence, when $n\rightarrow\infty$,
\[
A(\ell,n,k_n)\sim
\frac{\sqrt{\ell-1}}{2\pi}\times (n-1)!\times
\exp\left(
\mathscr{S}_{\ell}(n,k_n)+k_n\ \mathscr{E}_{\ell}\left(\frac{k_n}{n}\right)
\right)\ ,
\]
where
\[
\mathscr{E}_{\ell}(u):=H_1 u+H_2 u^2+\cdots+H_{\ell-1} u^{\ell-1}\ .
\]
\end{theorem}

Using bold and very interesting heuristics, Starr discovered the asymptotic formula in the above theorem and formulated it as Conjecture 2.3 in~\cite{Starr}. His formula, for any $\ell\ge 2$, allowed for logarithmic terms $u^a(-\ln u)^b$ in the correction $\mathscr{E}_{\ell}$ to the main quantity $\mathscr{S}_{\ell}$.
However, we found that no logarithms are present as soon as $\ell\ge 4$. The cases $\ell=2$ and $\ell=3$ are treated separately in the next two theorems, and they indeed involve logarithms. Altogether, Theorems \ref{mainthm1}, \ref{mainthm2}, and \ref{mainthm3} thus provide a complete proof for Starr's conjecture.
The case of pairs of permutations is the object of the next main result.

\begin{theorem}\label{mainthm2}
There exist two explicit constants $I_{1,{\rm log}},I_{1}$, such that for all $s\in(0,1)$, and all sequences $(k_n)_{n\ge 1}$, with
$1\le k_n\le n$ that are such that
$k_n\sim s n^{\frac{1}{2}}$, we have the asymptotic equivalence, when $n\rightarrow\infty$,
\[
A(2,n,k_n)\sim
\frac{1}{2\pi}\times (n-1)!\times
\exp\left(
\mathscr{S}_{2}(n,k_n)+k_n\ \mathscr{E}_{2}\left(\frac{k_n}{n}\right)
\right)\ ,
\]
where
\[
\mathscr{E}_{2}(u):=I_{1,{\rm log}}\ u \ln u+I_1 u\ .
\]
\end{theorem}

The case of triples of permutations is the object of our third main result.

\begin{theorem}\label{mainthm3}
There exist three explicit constants $J_1, J_{2,{\rm log}},J_{2}$, such that for all $s\in(0,1)$, and all sequences $(k_n)_{n\ge 1}$, with
$1\le k_n\le n$ that are such that
$k_n\sim s n^{\frac{2}{3}}$, we have the asymptotic equivalence, when $n\rightarrow\infty$,
\[
A(3,n,k_n)\sim
\frac{1}{\pi\sqrt{2}}\times (n-1)!\times
\exp\left(
\mathscr{S}_{3}(n,k_n)+k_n\ \mathscr{E}_{3}\left(\frac{k_n}{n}\right)
\right)\ ,
\]
where
\[
\mathscr{E}_{3}(u):=J_1 u
+J_{2,{\rm log}}\ u^2 \ln u+J_2 u\ .
\]
\end{theorem}

Since some are quite involved, the definitions of the $H,I,J$ constants are deferred to the next section \S\ref{constdefsec}. Note that our use of the notation $\sim$, $o(\cdot)$, and $O(\cdot)$ conforms to standard practice as in~\cite{deBruijn}. In particular, our statements of asymptotic equivalence $f(p)\sim g(p)$, when a parameter $p$ goes to some limit, mean that the limit of the ratio $\frac{f(p)}{g(p)}$ is $1$. We will also use $\sim$ for full asymptotic series expansions.

As a corollary of the previous theorems, we obtain strict log-concavity with respect to $k$ of the $A(\ell,n,k)$ in the hitherto unexplored regime of typical values of $k$, in the light of the central limit theorem proved in~\cite{AbdesselamS}.

\begin{corollary}\label{maincor}
For any $\ell\ge 2$, any $s\in(0,1)$, and all sequences $(k_n)_{n\ge 1}$, with
$1\le k_n\le n$ that are such that
$k_n\sim s n^{\frac{\ell-1}{\ell}}$, we have
\[
A(\ell,n,k_n)^2>A(\ell,n,k_n-1)A(\ell,n,k_n+1)\ ,
\]
for $n$ large enough.
\end{corollary}

As oulined in \S\ref{outlinesec}, our proof of the main three theorems is primarily based on bivariate saddle point analysis. It also crucially relies on detailed asymptotics, when $t\rightarrow 0^{+}$, of some remarkable multiple sums already considered in ~\cite{AbdesselamS} and denoted by $Z_{\alpha_1,\ldots,\alpha_{\ell}}^{[\ell]}(t)$. 
Now letting $\ell$ be arbitrary in $\mathbb{N}$, for all $(\alpha_1,\ldots,\alpha_{\ell})\in\mathbb{C}^{\ell}$, and all $t\in (0,\infty)$, we define these quantities by
\begin{equation}
Z_{\alpha_1,\ldots,\alpha_{\ell}}^{[\ell]}(t):=\sum_{\delta_1,\ldots,\delta_{\ell}=1}^{\infty}
\delta_{1}^{\alpha_1-1}\cdots\delta_{\ell}^{\alpha_{\ell}-1}\ e^{-\delta_1\cdots\delta_{\ell} t}\ .
\label{Zdefeq}
\end{equation}
The degenerate case $Z_{\varnothing}^{[0]}=e^{-t}$ is not excluded from our discussion. The case $Z_{\alpha}^{[1]}(t)={\rm Li}_{1-\alpha}(e^{-t})$
is related to polylogarithms. Higher $\ell$ cases of these $Z$ functions correspond to  evaluations at $z=e^{-t}$ of ordinary generating functions of multiplicative functions of number theory given by multiple Dirichlet convolution of power laws. These occur in many applications, e.g., $Z_{1,1}^{[2]}$ is the Lambert series studied by Wigert~\cite{Wigert}, and used for proving a mean value theorem for the Riemann zeta function~\cite[\S7.15]{Titchmarsh}. The particular case $Z_{1,0}^{[2]}$ is intimately connected to the Hardy-Ramanujan asymptotics for the partition function $p(n)$.
Our last main result is a reasonably explicit full asymptotic expansion as $t\rightarrow 0^{+}$ of the $Z$ functions for all $\ell$ and all complex numbers $\alpha_1,\ldots,\alpha_{\ell}$, including situations with multiplicities or resonances with the nonpositive integers which account for logarithmic factors.
The next theorem will come as no surprise to analytic number theorists and experts in the Mellin approach to asymptotic expansions (see, e.g.,~\cite{FlajoletGD} and~\cite{Zagier}). We believe it can be useful to applied mathematicians, as well as theoretical and mathematical physicists, to have such asymptotics collected in one place.
By $[z^k] f(z)$ we will denote the coefficient of $z^k$ in the power series expansion of $f(z)$ at $z=0$. We will also use $\bbone\{\cdots\}$ for the indicator function of the enclosed logical statement. As usual, $\Gamma(z)$ will denote the Euler gamma function.

\begin{theorem}\label{mainthm4}
For any $\ell\ge 0$, and $(\alpha_1,\ldots,\alpha_{\ell})\in\mathbb{C}^{\ell}$,
as $t\rightarrow 0^{+}$, we have the full asymptotic expansion
\[
Z_{\alpha_1,\ldots,\alpha_{\ell}}^{[\ell]}(t)\sim\sum_{(a,b)}
\mathscr{C}_{a,b}^{[\ell]}(\alpha_1,\ldots,\alpha_{\ell})\ 
t^a(-\ln t)^b\ ,
\]
where the sum is over pairs where $a\in\mathbb{N}\cup\{-\alpha_1,-\alpha_2,\ldots,-\alpha_{\ell}\}$, and $b\in\{0,1,\ldots,\ell\}$, and where
\begin{eqnarray}
\mathscr{C}_{a,b}^{[\ell]}(\alpha_1,\ldots,\alpha_{\ell}) &:=&
\frac{1}{b!}\ \sum_{p,k_1,\ldots,k_{\ell}\ge 0}
\bbone\{b+p+k_1+\cdots+k_{\ell}=\ell\} 
\left([z^p]\ z\Gamma(-a+z)\right) \nonumber \\
 & & \qquad\qquad\qquad\times\prod_{j=1}^{\ell}\left(
[z^{k_j}]\ z\zeta(-a+1-\alpha_j+z)
\right)\ .
\label{genCeq}
\end{eqnarray}
\end{theorem}

The precise meaning of the asymptotic expansion above is that for all $N\ge 0$, the difference between the $Z$ function on the left-hand side and the sum on right-hand side restricted to pairs $(a,b)$ with ${\rm Re}\ a\le N$ is $o(t^N)$. This therefore is a complete asymptotic expansion to order $O(t^{\infty})$. As in~\cite{AbdesselamS}, we will use, for $\ell\ge 1$ and $m\in\mathbb{Z}$, the abbreviated notation $Z_{m}^{[\ell]}(t)$ for $Z$ functions with a descending staircase pattern. 
Namely, we let
\begin{equation}
Z_{m}^{[\ell]}(t):=Z_{m,m-1,\ldots,m-\ell+1}^{[\ell]}(t)\ .
\label{abbZdefeq}
\end{equation}
For the latter functions, we will provide much more explicit asymptotic expansion formulas in \S\ref{asympsec}. The most important ones for the proof of Theorems \ref{mainthm1}, \ref{mainthm2}, and \ref{mainthm3} are $Z_{\ell}^{[\ell]}$ and $Z_{\ell-1}^{[\ell]}$. We will also need leading asymptotics for $Z_{\ell+1}^{[\ell]}$ and $Z_{\ell+2}^{[\ell]}$ which were already given in~\cite[\S2]{AbdesselamS}.

\subsection{Definition of the constants featuring in the main theorems}
\label{constdefsec}

We first address the definition of the constants featuring in Theorem \ref{mainthm1}, when $\ell\ge 4$.
We define the sequence of constants $(A_j)_{j\ge 0}$ by letting
\[
A_0:=1\ ,\qquad A_1:=\frac{\zeta(0)}{(\ell-1)\zeta(\ell)}\ ,
\qquad A_2:=\frac{\zeta(-1)\zeta(0)}{(\ell-1)(\ell-2)\zeta(\ell-1)\zeta(\ell)}\ ,
\]
and then $A_j:=0$, for $j\ge 3$.
We also define the sequence of constants $(B_j)_{j\ge 0}$ by letting
\[
B_0:=1\ ,\qquad B_1:=\frac{\zeta(0)}{(\ell-2)\zeta(\ell)}\ ,
\qquad B_2:=\frac{\zeta(-1)\zeta(0)}{(\ell-2)(\ell-3)\zeta(\ell-1)\zeta(\ell)}\ ,
\]
as well as
\[
B_{\ell-1}:=\left\{
\begin{array}{cll}
\frac{\zeta'(-2)\zeta(-1)\zeta(0)}{2\zeta(2)\zeta(3)\zeta(4)} & {\rm if} & \ell=4 \ , \\
\frac{\zeta(-3)\zeta'(-2)\zeta(-1)\zeta(0)}{6\zeta(2)\zeta(3)\zeta(4)\zeta(5)} & {\rm if} & \ell=5\ , \\
0 & {\rm if} & \ell\ge 6\ ,
\end{array}
\right.
\]
and then $B_j:=0$, for $j\notin\{0,1,2,\ell-1\}$.

We next define the sequence of constants $(C_j)_{j\ge 1}$ by letting
\[
C_j:=\sum_{m\ge 0} B_m\sum_{k\ge 0}(-1)^k
\sum_{j_1,\ldots,j_k\ge 1}\bbone\{m+j_1+\cdots+j_k=j-1\}\ A_{j_1}\cdots A_{j_k}\ ,
\]
for all $j\ge 1$.

Then we introduce the sequence of constants $(D_p)_{p\ge 1}$ defined by letting, for all $p\ge 1$,
\[
D_p:=\frac{1}{p!}\sum_{\mathbf{m}}
\bbone\left\{\sum_{j\ge 2}(j-1)m_j=p-1\right\}
\times\frac{\left(\sum_{j\ge 2} j m_j\right)!}{\prod_{j\ge 2}m_j!}
\times (-1)^{\sum_{j\ge 1}m_j}
\times \prod_{j\ge 2} C_{j}^{m_j}\ ,
\] 
where the sum is over infinite yet almost finite sequences of nonnegative integers $\mathbf{m}=(m_j)_{j\ge 2}$, i.e., which satisfy $m_j=0$ for $j$ large enough.

We also define the sequence of constants $(E_j)_{j\ge 1}$ by letting
\[
E_j:=\sum_{k,j_1,\ldots,j_{\ell}\ge 0}
\bbone\{k+j_1+\cdots+j_{\ell}=j\}\ A_k C_{j_1+1}\cdots C_{j_{\ell}+1}\ ,
\]
for all $j\ge 1$.
This in turn allows us to define
\[
F_j:=(\ell-1)^j\times
\sum_{k\ge 1} E_k\sum_{j_1,\ldots,j_k\ge 1}
\bbone\{j_1+\cdots+j_k=j\}\ D_{j_1}\cdots D_{j_k}\ ,
\]
for all $j\ge 1$, as well as
\[
G_j:=
\sum_{k\ge 1} \frac{(-1)^{k-1}}{k}
\sum_{j_1,\ldots,j_k\ge 1}
\bbone\{j_1+\cdots+j_k=j\}\ F_{j_1}\cdots F_{j_k}\ ,
\]
also for all $j\ge 1$.
Finally, we can let, for all $j\ge 1$,
\[
H_j:=(\ell-1)^{j+1}\ D_{j+1}+G_j\ ,
\]
although only the first $\ell-1$ values are needed in the statement of Theorem \ref{mainthm1}.

For $\ell=2$, the needed constants in the statement of Theorem \ref{mainthm2} are
\[
I_{1,{\rm log}}:=-\frac{\zeta(0)}{\zeta(2)}=\frac{3}{\pi^2}\ ,
\]
and
\[
I_{1}:=\frac{\zeta'(0)}{\zeta(2)}=-\frac{3\ln(2\pi)}{\pi^2}\ ,
\]
where we used known special values of the zeta function and its derivative (see, e.g.,~\cite[Thm. 1.2.4, Cor. 1.3.2, and Cor. 1.3.3]{AndrewsAR}).

For $\ell=3$, the needed constants in the statement of Theorem \ref{mainthm3} are
\[
J_1:=\frac{2\zeta(0)}{\zeta(3)}=-\frac{1}{\zeta(3)}\ ,
\]
\[
J_{2,{\rm log}}:=-\frac{4\zeta(-1)\zeta(0)}{\zeta(2)\zeta(3)}
=-\frac{1}{\pi^2\zeta(3)}\ ,
\]
and
\[
J_2:=\frac{1}{\zeta(2)\zeta(3)^2}\left[
-3\zeta(0)^2\zeta(2)-2(\ln 2)\zeta(-1)\zeta(0)\zeta(3)
+4\zeta'(-1)\zeta(0)\zeta(3)+4\zeta(-1)\zeta'(0)\zeta(3)
\right]\ .
\]
In the last rather formidable constant, we left the zetas unevaluated since, not only $\zeta(3)$ has no simpler expression, but neither does $\zeta'(-1)$ which is related to the Glaisher-Kinkelin constant.

\subsection{Outline of our approach}
\label{outlinesec}

For any $\ell\ge 1$, the numbers $A(\ell,n,k)$ can be packaged in the bivariate generating function
\[
\mathscr{G}_{\ell}(x,z):=\sum_{n=0}^{\infty}\sum_{k=0}^{n}
\frac{1}{n!} A(\ell,n,k) x^{k}z^{n}\ .
\]
From a result by Bryan and Fulman~\cite{BryanF} (see also~\cite{AbdesselamBDV}), we have the explicit infinite product representation
\begin{equation}
\mathscr{G}_{\ell}(x,z)=\prod_{\delta_1,\ldots,\delta_{\ell-1}=1}^{\infty}
\left(1-z^{\delta_1\cdots\delta_{\ell-1}}\right)^{-x\delta_1^{\ell-2}\delta_2^{\ell-3}\cdots\delta_{\ell-2}}\ .
\label{BFeq}
\end{equation}
This formula still makes sense and is correct if $\ell=1$, provided one interprets empty products as equal to $1$, and empty tuples as the function $\delta:\varnothing\rightarrow\mathbb{Z}_{>0}$ with empty graph. In this case, the above reads
\[
\sum_{n=0}^{\infty}\frac{x(x+1)\cdots(x+n-1)}{n!}\ z^n=(1-z)^{-x}\ ,
\]
which agrees with the definition of the $A(1,n,k)$ as the unsigned Stirling numbers of the first kind $|s(n,k)|$. The proof of (\ref{BFeq}) in~\cite{BryanF} was in the sense of formal power series in $\mathbb{C}[[x,z]]$, however, as explained in~\cite[\S3]{AbdesselamS}, the equality holds at the level of analytic functions. Indeed, using the same notations as~\cite{AbdesselamS}, let us define for $|z|<1$,
\begin{eqnarray*}
\mathscr{L}_{\ell}(z) & := & -\sum_{\delta_1,\ldots,\delta_{\ell-1}=1}^{\infty}
\delta_1^{\ell-2}\delta_2^{\ell-3}\cdots\delta_{\ell-2}
\ {\rm Log}(1-z^{\delta_1\cdots\delta_{\ell-1}}) \\
 & = & \sum_{\delta_1,\ldots,\delta_{\ell}=1}^{\infty}
\delta_1^{\ell-2}\delta_2^{\ell-3}\cdots\delta_{\ell-2}\delta_{\ell-1}^{0}\delta_{\ell}^{-1}
\ z^{\delta_1\cdots\delta_{\ell}}\ ,
\end{eqnarray*}
where $\rm Log$ is the principal branch of the complex logarithm. The two series representations for $\mathscr{L}_{\ell}(z)$ converge absolutely when $|z|<1$, and therefore the product in (\ref{BFeq}) is an absolutely convergent one, and defines a jointly analytic function of $(x,z)$ in the domain $\mathbb{C}\times\{z\in\mathbb{C}\ |\ |z|<1\}$ (only the case of real $x$ was considered in~\cite[\S3]{AbdesselamS}, but the same works for any complex value of $x$). Hence, $\mathscr{G}_{\ell}(x,z)$ is analytic on the same domain and can be rewritten
\[
\mathscr{G}_{\ell}(x,z)=\exp\left(x\mathscr{L}_{\ell}(z)\right)\ .
\]
As a result, we can extract our coefficients of interest using two Cauchy integrals
\[
\frac{A(\ell,n,k)}{n!}=\oint_{|z|=e^{-t}}\frac{{\rm d}z}{2i\pi z}
\ \oint_{|x|=\rho}\frac{{\rm d}x}{2i\pi x}
\ \ z^{-n}x^{-k}\ \mathscr{G}_{\ell}(x,z)\ ,
\]
with counterclockwise contours given by circles around the origin with radius $e^{-t}$ for the $z$ variable, and radius $\rho$ for the $x$ variable. For the moment, $t$ and $\rho$ are arbitrary in $(0,\infty)$, but they will have to both be optimized in order to study the above double integral using bivariate saddle point analysis.
Apart from dealing with a double instead of a single integral, we follow the general philosophy explained, e.g., in~\cite[\S5.1]{deBruijn}, namely we will pick the optimal two-dimensional contour in order to minimize the maximum over the contour of the modulus of the integrand. Writing $z=e^{-t+i\theta}$ and $x=\rho e^{i\omega}$, with $\theta,\omega\in(-\pi,\pi)$, we consider
\[
\mathscr{M}_{\ell}(t,\rho):=\sup_{\theta,\omega\in(-\pi,\pi)}
\left|z^{-n}x^{-k}\ \mathscr{G}_{\ell}(x,z)\right|\ .
\]
We note that
\[
\ln \left|z^{-n}x^{-k}\ \mathscr{G}_{\ell}(x,z)\right|
=nt-k\ln\rho+{\rm Re}\left(x\mathscr{L}_{\ell}(z)\right)\ ,
\]
and
\[
{\rm Re}\left(x\mathscr{L}_{\ell}(z)\right)\le 
\left|x\mathscr{L}_{\ell}(z)\right|
\le \rho
\sum_{\delta_1,\ldots,\delta_{\ell}=1}^{\infty}
\delta_1^{\ell-2}\delta_2^{\ell-3}\cdots\delta_{\ell-2}\delta_{\ell-1}^{0}\delta_{\ell}^{-1}
\ |z|^{\delta_1\cdots\delta_{\ell}}\ .
\]
We thus see that the maximum is reached at $\theta=\omega=0$, and that
\begin{equation}
\ln\mathscr{M}_{\ell}(t,\rho)=nt-k\ln \rho+\rho Z_{\ell-1}^{[\ell]}(t)\ ,
\label{logMeq}
\end{equation}
where we took advantage of the notation introduced in (\ref{Zdefeq}) and (\ref{abbZdefeq}).
The optimal $t$ and $\rho$ are obtained by solving the system
\[
\left\{
\begin{array}{ccc}
\frac{\partial}{\partial t}\ln \mathscr{M}_{\ell}(t,\rho) & = & 0\ , \\
\frac{\partial}{\partial\rho}\ln \mathscr{M}_{\ell}(t,\rho) & = & 0\ ,
\end{array}
\right.
\]
i.e.,
\[
\left\{
\begin{array}{ccc}
n & = & \rho Z_{\ell}^{[\ell]}(t)\ , \\
k & = & \rho Z_{\ell-1}^{[\ell]}(t)\ .
\end{array}
\right.
\]
We refer to~\cite[\S2]{AbdesselamS} for the basic properties of the $Z$ functions, and in particular their strict positivity and smoothness on $(0,\infty)$, and the relation $\frac{{\rm d}}{{\rm d}t}Z_{\ell-1}^{[\ell]}(t)=-Z_{\ell}^{[\ell]}(t)$ used above.
We define the function
\begin{equation}
h_{\ell}(t):=(\ell-1)\times\frac{Z_{\ell-1}^{[\ell]}(t)}{Z_{\ell}^{[\ell]}(t)}\ ,
\label{hdefeq}
\end{equation}
where the $\ell-1$ factor is introduced as a matter of convenience, so that  $h_{\ell}(t)\sim t$ when $t\rightarrow 0^{+}$. We will show in \S\ref{asympsec} that $h_{\ell}$ is a bijection from $(0,\infty)$ onto $(0,\ell-1)$ and therefore has a well defined inverse function $h_{\ell}^{-1}:(0,\ell-1)\rightarrow(0,\infty)$.

Now consider a sequence $(k_n)_{n\ge 1}$ as in the statement of Theorems \ref{mainthm1}, \ref{mainthm2}, and \ref{mainthm3}.  For $n$ large enough, we have $0<k<n$ because $s\in(0,1)$.
Simple algebra shows that the unique solution $t_n,\rho_n$ to the system
\[
\left\{
\begin{array}{ccc}
n & = & \rho_n Z_{\ell}^{[\ell]}(t_n)\ , \\
k_n & = & \rho_n Z_{\ell-1}^{[\ell]}(t_n)\ ,
\end{array}
\right.
\]
is obtained by letting
\begin{equation}
t_n:=h_{\ell}^{-1}\left(\frac{(\ell-1)k_n}{n}\right)\ ,
\label{tndefeq}
\end{equation}
and then
\begin{equation}
\rho_n:=\frac{n}{Z_{\ell}^{[\ell]}(t_n)}\ ,
\label{rhondefeq}
\end{equation}
which we will assume for the rest of the proof of the main theorems.
The above choice of contours gives us the multiplicative decomposition
\begin{equation}
\frac{A(\ell,n,k_n)}{n!}=\mathscr{M}_{\ell}(t_n,\rho_n)\times
\mathscr{I}_n\ ,
\label{multdeceq}
\end{equation}
where the second factor is the integral
\begin{eqnarray}
\mathscr{I}_n &:= & \oint_{|z|=e^{-t_n}}\frac{{\rm d}z}{2i\pi z}
\ \oint_{|x|=\rho_n}\frac{{\rm d}x}{2i\pi x}
\ \ \frac{z^{-n}x^{-k_n}\ \mathscr{G}_{\ell}(x,z)}
{e^{n t_n}\rho_n^{-k_n}\ \mathscr{G}_{\ell}(\rho_n,e^{-t_n})} \label{curlyIdefeq}\\
 &= & \int_{-\pi}^{\pi}\frac{{\rm d}\theta}{2\pi}
\ \int_{-\pi}^{\pi}\frac{{\rm d}\omega}{2\pi}
\ \ \mathcal{I}_n(\theta,\omega)\ ,
\nonumber
\end{eqnarray}
with the last integrand given by
\begin{equation}
\mathcal{I}_n(\theta,\omega) :=
\frac{e^{-in\theta} e^{-ik_n\omega}
\exp\left(\rho_n e^{i\omega}\mathscr{L}_{\ell}
\left(e^{-t_n+i\theta}\right)\right)}{\exp\left(\rho_n \mathscr{L}_{\ell}
\left(e^{-t_n}\right)\right)}\ .
\label{Indefeq}
\end{equation}
By design, the maximum of $|\mathcal{I}_n(\theta,\omega)|$ is equal to $1$ and is achieved at $\theta=\omega=0$.
Our proof of the main theorems rely on precise asymptotics for the prefactor $\mathscr{M}_{\ell}(t_n,\rho_n)$ using the results of \S\ref{asympsec}, and careful bivariate saddle point analysis for the integral $\mathscr{I}_n$ provided in \S\ref{saddlesec}.

\subsection{Related work}

The main motivation for the present article is the conjecture by the author~\cite{AbdesselamAC} which says that the $A(\ell,n,k)$ should be log-concave with respect to $k$. In the particular case $\ell=2$, this conjecture is older and due to Heim and Neuhauser~\cite[Challenge 3]{HeimN}.
For $\ell=2$, as explained in~\cite{AbdesselamBDV}, the equation (\ref{BFeq}) reduces to
\[
\sum_{n=0}^{\infty}\sum_{k=0}^{n}
\frac{1}{n!} A(2,n,k) x^{k}z^{n}=
\prod_{\delta=1}^{\infty}
(1-z^{\delta})^{-x}\ ,
\]
which relates to the D'Arcais polynomials and the Nekrasov-Okounkov formula.
We believe Corollary \ref{maincor} is the first result which establishes the log-concavity for typical values of $k$.

A central problem in analytic number theory is to obtain estimates on partial sums $\sum_{n\le x}f(n)$ where $f(n)$ is a function of number-theoretic interest. 
Strict partial sum with a sharp cutoff are difficult to analyze. One typically has to modify the definition to $\sum'$ so the last term counts for half if $x$ is an integer.
Such partial sums can also be written $\sum_{n}f(n)\phi(tn)$ where $t=\frac{1}{x}$ and $\phi(t)=\bbone\{0<t\le 1\}$. Smooth summatory functions like $\phi(t)=e^{-t}$ are much nicer and more amenable to the derivation of full asymptotic expansions as in Theorem \ref{mainthm4} (see the insightful discussion in~\cite[\S3.7]{Tao} and~\cite[p. 33]{Kowalski}). The above sums were also considered, as some renormalization of an a priori divergent series $\sum_{n=1}^{\infty}f(n)$, by Ramanujan (see~\cite{BerndtE}). The two main techiniques used for deriving asymptotics for these sums are the real variable methods based on the Euler-MacLaurin formula, and the complex variable methods based on the Mellin transform. One could argue that the switch from the former to the latter, as a preferred choice, can be traced back to Riemann's article on the zeta function~\cite{Riemann} (see also the appendix of~\cite{Edwards} for an English translation).
Although the proofs of Theorem \ref{mainthm4} and the more explict Proposition \ref{Zmprop} follow the standard Mellin transform approach~\cite{FlajoletGD,Zagier}, we believe they could be useful, given that such asymptotics are needed for many applications (see, e.g.,~\cite{NinhamHFG}).

The main method we use in order to prove Theorems \ref{mainthm1}, \ref{mainthm2} and \ref{mainthm3}, is saddle point analysis. We follow similar steps to those in~\cite{TenenbaumWL,DebruyneT}, for example. However, here we need bivariate saddle point analysis. We are in a lucky situation where the two-dimensional contour of integration, which is a torus, does not undergo a change of topology. All we have to do is adjust the radius of each of the two circles involved, and we can thus avoid the more high-powered techniques, like Morse theory, sometimes required for multivariate saddle point problems~\cite{PemantleW}.
For $\ell=1$, i.e., Stirling numbers of the first kind, asymptotics have been studied by many authors, see in particular the work of Moser and Wyman~\cite{MoserW}.
However, for $\ell\ge 2$, we are not aware of rigorous asymptotics for the $A(\ell,n,k)$ before the present article.
Our bivariate saddle point analysis shares some thematic similarity to the works~\cite{Hwang,LipnikMT2}. We note that the symmetrization trick we used for proving the crucial log-convexity lemma below appeared in~\cite{MadritschW,LipnikMT1} (see also~\cite{AbdesselamCorr}).

\section{Complete asymptotics and log-convexity for the $Z$ functions}
\label{asympsec}

\noindent{\bf Proof of Theorem \ref{mainthm4}:}
We start with the Mellin inversion formula (see~\cite[Cor 0.14, p. 177]{TenenbaumBook}),
\[
e^{-u}=\frac{1}{2i\pi}\int_{\sigma_{\rm ini}-i\infty}^{\sigma_{\rm ini}+i\infty}
\Gamma(w) u^{-w}\ {\rm d}w
\]
valid for $u>0$ and $\sigma_{\rm ini}>0$. We use it for $u=\delta_1\cdots\delta_{\ell}t$ and insert this in the definition of the $Z$ functions. Namely, we write for any fixed $t>0$,
\[
Z_{\alpha_1,\ldots,\alpha_{\ell}}^{[\ell]}(t)=\sum_{\delta_1,\ldots,\delta_{\ell}=1}^{\infty}
\delta_{1}^{\alpha_1-1}\cdots\delta_{\ell}^{\alpha_{\ell}-1}\times 
\frac{1}{2i\pi}\int_{\sigma_{\rm ini}-i\infty}^{\sigma_{\rm ini}+i\infty}
\Gamma(w) (\delta_1\cdots\delta_{\ell} t)^{-w}\ {\rm d}w\ .
\]
By the complex Stirling asymptotic formula for the Gamma function in the vertical direction (see~\cite[Cor 1.4.4]{AndrewsAR} or~\cite[Cor 0.13, p. 176]{TenenbaumBook}), one has a bound $|\Gamma(w)|\le C e^{-c|{\rm Im} w|}$, 
for $|{\rm Im}\ w|ge 1$,
with constants $C,c>0$ depending on the abscissa $\sigma_{\rm ini}$ of the initial contour of integration.
Since
\[
\left|\Gamma(w)\ \delta_{1}^{\alpha_1-1-w}\cdots\delta_{\ell}^{\alpha_{\ell}-1-w} \ t^{-w} 
\right|\le C\ \delta_{1}^{({\rm Re}\ \alpha_1)-1-\sigma_{\rm ini}}\cdots
\delta_{\ell}^{({\rm Re}\ \alpha_{\ell})-1-\sigma_{\rm ini}}
\ t^{-\sigma_{\rm ini}}\ e^{-c|{\rm Im} w|}\ ,
\]
we get a finite result when integrating over $w$ and summing over the $\delta$'s, provided we made sure to pick $\sigma_{\rm ini}>\max(0, {\rm Re}\ \alpha_1,\ldots,{\rm Re}\ \alpha_{\ell})$, in order for all series to converge.
Then, Fubini's theorem applies and gives, after commuting the sums and the integral,
\begin{eqnarray*}
Z_{\alpha_1,\ldots,\alpha_{\ell}}^{[\ell]}(t) & = & \frac{1}{2i\pi}
\int_{\sigma_{\rm ini}-i\infty}^{\sigma_{\rm ini}+i\infty}
\sum_{\delta_1,\ldots,\delta_{\ell}=1}^{\infty}
\Gamma(w)
\ \delta_{1}^{\alpha_1-1-w}\cdots\delta_{\ell}^{\alpha_{\ell}-1-w}
\ t^{-w}\ {\rm d}w \\
 & = & \frac{1}{2i\pi}
\int_{\sigma_{\rm ini}-i\infty}^{\sigma_{\rm ini}+i\infty}
\Gamma(w)\ \zeta(w+1-\alpha_1)\cdots\zeta(w+1-\alpha_{\ell})\ t^{-w}\ {\rm d}w \ .
\end{eqnarray*}
It is known that the Gamma function $\Gamma(w)$ is a meromorphic function on the complex plane with simple poles at $w\in\mathbb{Z}_{\le 0}:=\{0,-1,-2,\ldots\}$. The Riemann zeta function also admits analytic continuation to $\mathbb{C}$ except for a simple pole at $w=1$. We refer to~\cite[Ch. 1]{AndrewsAR} for these standard facts. In the case of the Riemann zeta function, a lot of specific information can also be found in the books~\cite{Apostol,Edwards,TenenbaumBook,Titchmarsh}.
The function $w\mapsto t^{-w}=e^{(-\ln t)w}$ is entire analytic.
The above integrand therefore is meromorphic with potential poles contained in the discrete subset of $\mathbb{C}$ given by $\mathcal{P}:=\mathbb{Z}_{\le 0}\cup\{\alpha_1,\ldots,\alpha_{\ell}\}$. We said potential poles, because these could be cancelled by the mysterious zeros of the zeta function, depending on the placement of the points $\alpha_1,\ldots,\alpha_{\ell}$ in the complex plane.
We can now move the contour, to the left, to a final abscissa $\sigma_{\rm fin}<\sigma_{\rm ini}$, while collecting the contributions of the poles crossed via Cauchy's residue theorem.
We first replace the improper integral $\int_{\sigma_{\rm ini}-i\infty}^{\sigma_{\rm ini}+i\infty}$ by $\lim_{T\rightarrow\infty}\int_{\sigma_{\rm ini}-iT}^{\sigma_{\rm ini}+iT}$. We apply the residue theorem to the counterclockwise rectangular contour going successively through the points $\sigma_{\rm ini}-iT$, $\sigma_{\rm ini}+iT$, $\sigma_{\rm fin}+iT$, $\sigma_{\rm fin}-iT$ and then back to $\sigma_{\rm ini}-iT$.
The contribution of the two segments with ${\rm Im}\ w=\pm T$ go to zero, because the previous exponential decay in the vertical direction for $|\Gamma(w)|$ can be made uniform in strips $\sigma_{\rm fin}\le {\rm Re}\ w\le \sigma_{\rm ini}$ with finite width. Moreover, on such strips the zeta function grows uniformly at most polynomially, i.e., $\exists K,\kappa>0$, such that for all $w$ with $\sigma_{\rm fin}\le {\rm Re}\ w\le \sigma_{\rm ini}$ and $|{\rm Im}\ w|\ge 1$, we have
$|\zeta(w)|\le K |{\rm Im}\ w|^{\kappa}$ (see~\cite[Thm 12.23]{Apostol}
or~\cite[\S3.4]{TenenbaumBook}).
Trivially, one can deduce similar uniform polynomial growth bounds along strips of finite
width for products of shifted zeta functions as features in the above integral representation for $Z$.
As a result, provided the negative real number $\sigma_{\rm fin}$ does not belong to the set $\mathbb{Z}_{\le 0}\cup\{{\rm Re}\ \alpha_1,\ldots,{\rm Re}\ \alpha_{\ell}\}$, the moving of the contour is justified and yields
\begin{eqnarray*}
Z_{\alpha_1,\ldots,\alpha_{\ell}}^{[\ell]}(t) & = & 
\sum_{\substack{
q\in\mathcal{P} \\
\sigma_{\rm fin}<{\rm Re}\ q<\sigma_{\rm ini}}}
{\rm res}_{w=q}\left(\ 
\Gamma(w)\ \zeta(w+1-\alpha_1)\cdots\zeta(w+1-\alpha_{\ell})\ e^{(-\ln t)w}
\ \right) \\
 & & +\frac{1}{2i\pi}
\int_{\sigma_{\rm fin}-i\infty}^{\sigma_{\rm fin}+i\infty}
\Gamma(w)\ \zeta(w+1-\alpha_1)\cdots\zeta(w+1-\alpha_{\ell})\ t^{-w}\ {\rm d}w \ .
\end{eqnarray*}
By the bounds discussed above, the last integral remainder term is absolutely convergent and is bounded by a constant times $t^{-\sigma_{\rm fin}}$ which is $O(t^N)$ if $\sigma_{\rm fin}<-N$.
To compute the residue at $q$ we change variables to $z=w-q$ and we multiply the function by $z^{\ell+1}$. This converts the operation of taking the residue into that of taking the coefficient of $z^{\ell}$, i.e., computing
\[
[z^{\ell}]\left(
\ z\Gamma(q+z)\ z\zeta(q+z+1-\alpha_1)\cdots z\zeta(q+z+1-\alpha_{\ell}) 
\ t^{-q}\ e^{(-\ln t)z}\ \right)\ .
\]
Expanding the exponential as $\sum_{k=0}^{\infty}\frac{(-\ln t)^k}{k!}$, the function $z\Gamma(q+z)$ as a sum over $p\ge 0$, and the factors $z\zeta(q+z+1-\alpha_{j}) $
as sums over $k_j$, with $1\le j\le k$, completes the proof of Theorem \ref{mainthm4}.
\qed

We believe the multiplicative shift by factors of $z$ makes the computation of the residues easier, since one is dealing with convergent power series instead of Laurent series. It also makes the range of the summation indices $p,k_1,k_2,\ldots$ faster to determine. More explicit formulas for these residues can readily be obtained from the following expansions.

Let $u\in\mathbb{C}$. When $u\neq 1$, we immediately get by Taylor expansion, that for all $k\ge 0$,
\begin{equation}
[z^k]\ z\zeta(u+z)=\bbone\{k\ge 1\}\times\frac{\zeta^{(k-1)}(u)}{(k-1)!}\ .
\label{zetaregeq}
\end{equation}
For $u=1$, we have that for all $k\ge 0$,
\begin{equation}
[z^k]\ z\zeta(1+z)=\ \bbone\{k=0\}+
\bbone\{k\ge 1\}\times\frac{(-1)^{k-1}}{(k-1)!} \gamma_{k-1}\ ,
\label{zetasingeq}
\end{equation}
where $\gamma_0,\gamma_1,\gamma_2,\ldots$ are the Stieltjes constants (see~\cite[Thm. 3.2.1]{Lagarias}). In particular, $\gamma_0=\gamma$, the Euler-Mascheroni constant.

With regards to the gamma function, for $u\in\mathbb{C}\backslash\mathbb{Z}_{\le 0}$, we get by Taylor expansion that for all $k\ge 0$,
\begin{equation}
[z^k]\ z\Gamma(u+z)=\ \bbone\{k\ge 1\}\times\frac{\Gamma^{(k-1)}(u)}{(k-1)!}\ .
\label{gammaregeq}
\end{equation}
For $u=-a$, with $a\in\mathbb{N}$, we use the functional equation for $\Gamma$ to write
\begin{eqnarray*}
z\Gamma(u+z) & = & \frac{\Gamma(1+z)}{(u+z)(u+z+1)\cdots(-1+z)} \\
 & = & \frac{(-1)^a}{a!}
\times\frac{\Gamma(1+z)}{\left(1-z\right)\left(1-\frac{z}{2}\right)
\cdots\left(1-\frac{z}{a}\right)}\ ,
\end{eqnarray*}
from which we easily get, for all $k\ge 0$,
\begin{equation}
[z^k]\ z\Gamma(u+z)=\frac{(-1)^a}{a!}\times\sum_{m,k_1,\ldots,k_a\ge 0}
\bbone\{m+k_1+\cdots+k_a=k\}\times
\frac{\Gamma^{(m)}(1)}{m!\ 1^{k_1} 2^{k_2}\cdots a^{k_a}}\ .
\label{gammasingeq}
\end{equation}
We will put the above formulas to good use by deriving asymptotics, for the special case $Z_{m}^{[\ell]}$, $m\in\mathbb{Z}$, that are more explicit than what Theorem \ref{mainthm4} delivers.
For this we will need notation for the set
\[
S_{\ell,m}:=\{-m,-m+1,\ldots,-m+\ell-1\}\ .
\] 
\begin{proposition}\label{Zmprop}
For $m\in\mathbb{Z}$, and as $t\rightarrow 0^{+}$, we have the complete asymptotic expansion
\[
Z_{m}^{[\ell]}\sim\sum_{(a,b)\in (\mathbb{N}\cup S_{\ell,m})\times\{0,1\}}\mathscr{C}_{a,b}^{[\ell]}(m,m-1,\ldots,m-\ell+1)\ t^a 
(-\ln t)^b\ ,
\]
where the coefficients are given by the following formulas.

If $a\in\mathbb{N}\backslash S_{\ell,m}$, then
\[
\mathscr{C}_{a,b}^{[\ell]}(m,m-1,\ldots,m-\ell+1)=
\ \bbone\{b=0\}\times\frac{(-1)^a}{a!}\times
\prod_{j=1}^{\ell}\zeta(-a-m+j)\ .
\]
If $a\in S_{\ell,m}\backslash \mathbb{N}$, then
\[
\mathscr{C}_{a,b}^{[\ell]}(m,m-1,\ldots,m-\ell+1)=
\ \bbone\{b=0\}\times(-a-1)!\times
\prod_{\substack{j=1 \\ j\neq a+m+1}}^{\ell}\zeta(-a-m+j)\ .
\]
If $a\in\mathbb{N}\cup S_{\ell,m}$, then
\[
\mathscr{C}_{a,1}^{[\ell]}(m,m-1,\ldots,m-\ell+1)=
\ \frac{(-1)^a}{a!}\times \prod_{\substack{j=1 \\ j\neq a+m+1}}^{\ell}\zeta(-a-m+j)\ ,
\]
and
\begin{eqnarray*}
\mathscr{C}_{a,0}^{[\ell]}(m,m-1,\ldots,m-\ell+1) &=&
\frac{(-1)^a}{a!}\times \left(1+\frac{1}{2}+\cdots+\frac{1}{a}\right)\times
\prod_{\substack{j=1 \\ j\neq a+m+1}}^{\ell}\zeta(-a-m+j) \\
 & & +\sum_{\substack{j=1 \\ j\neq a+m+1}}^{\ell}
\zeta'(-a-m+j)\times \prod_{\substack{r=1 \\ r\notin\{a+m+1,j\}}}^{\ell}
\zeta(-a-m+r)\ .
\end{eqnarray*}
\end{proposition}

\noindent{\bf Proof:}

\noindent\underline{1st case:} Assume $a\in\mathbb{N}\backslash S_{\ell,m}$.
Then the zeta factors in (\ref{genCeq}) are regular at $z=0$ and the nonvanishing of the coefficient $[z^{k_j}]$ forces $k_1,\ldots,k_{\ell}\ge 1$. Nonnegativity and the constraint
$b+p+k_1+\cdots+k_{\ell}=\ell$ imply that the only contributing possibility is when $b=p=0$ and $k_1=\cdots=k_{\ell}=1$, i.e.,
\[
\mathscr{C}_{a,b}^{[\ell]}(m,m-1,\ldots,m-\ell+1) =
\bbone\{b=0\}\times
([z^0] z\Gamma(-a+z))\times\prod_{j=1}^{\ell}\left(
[z^1] z\zeta(-a-m+j+z)
\right)\ ,
\]
and the result follows from the use of formulas (\ref{gammasingeq}) and (\ref{zetaregeq}).

\noindent\underline{2nd case:} Assume
$a\in S_{\ell,m}\backslash \mathbb{N}$. This implies $a\in\mathbb{Z}_{<0}$
and $-m\le a\le -m+\ell-1$. The latter inequalities can be rewritten as
$1\le a+m+1\le\ell$, i.e., the value $a+m+1$ of $j$ which would lead to the pole of the zeta function is included in the range for $j$ in the product from (\ref{genCeq}).
If $j\in\{1,2,\ldots,\ell\}\backslash\{a+m+1\}$, the corresponding zeta factor is not evaluated at a pole which forces $k_j\ge 1$. Since $a<0$, we have that $-a$ is not a pole of the gamma function which forces $p\ge 1$. Hence
\[
p+\sum_{\substack{j=1 \\ j\neq a+m+1}}^{\ell}k_j\ge \ell
\]
which, given the various constraints at hand, leaves only one possibility, namely, $p=1$, $k_j=1$ for $j\neq a+m+1$, and $k_{a+m+1}=0$. Note that we must also have $b=0$, i.e., no logarithm. Again, using formulas
(\ref{gammaregeq}), (\ref{zetaregeq}), and (\ref{zetasingeq}),
we arrive at the result stated in the proposition.

\noindent\underline{3rd case:} Assume $a\in\mathbb{N}\cap S_{\ell,m}$.
Now the gamma function is at a pole, and among the zeta factors, exactly one corresponds to a pole, namely, the factor with $j=a+m+1$.
We must have $k_j\ge 1$ for $j\neq a+m+1$, which results in the cap
\[
b+p+k_{a+m+1}\le 1\ .
\]
This forces $b=0$
 or $b=1$, i.e., the staircase pattern with no multiplicities among the $\alpha$'s limits the power of the logarithm to at most 1, instead of the worst-case bound $\ell$ anticipated in Theorem \ref{mainthm4}.

If $b=1$, the only possibility for the summation indices is $p=k_{a+m+1}=0$, and $k_j=1$ for the other $j$'s which, using
(\ref{gammasingeq}), (\ref{zetaregeq}), and (\ref{zetasingeq}),
results in the evaluation
\[
\frac{(-1)^a}{a!}\times
\prod_{\substack{j=1 \\ j\neq a+m+1}}^{\ell}\zeta(-a-m+j)
\]
for the coefficient $\mathscr{C}_{a,1}^{[\ell]}$.

Let us now focus on the last possibility with $b=0$, where the pair $(p,k_{a+m+1})$ can be $(0,1)$, $(1,0)$, or $(0,0)$.
If $(p,k_{a+m+1})=(0,1)$, then we must also have $k_j=1$ for $j\neq a+m+1$, and
(\ref{gammaregeq}), (\ref{zetaregeq}), and (\ref{zetasingeq})
lead to the contribution
\[
\frac{(-1)^a}{a!}
\times\gamma\times
\prod_{\substack{j=1 \\ j\neq a+m+1}}^{\ell}\zeta(-a-m+j)
\]
to the coefficient $\mathscr{C}_{a,0}^{[\ell]}$.
If $(p,k_{a+m+1})=(0,1)$, then we must also have $k_j=1$ for $j\neq a+m+1$, and
(\ref{gammaregeq}), (\ref{zetaregeq}), and (\ref{zetasingeq})
lead to the contribution
\[
\frac{(-1)^a}{a!}
\times
\left(1+\frac{1}{2}+\cdots+\frac{1}{a}+\Gamma'(1)\right)
\times
\prod_{\substack{j=1 \\ j\neq a+m+1}}^{\ell}\zeta(-a-m+j)\ .
\]
If $(p,k_{a+m+1})=(0,0)$, then for exactly one value $r$ of $j\neq a+m+1$ we must have $k_r=1$ while the other $k_j$'s are zero. Moreover, we must sum over all possibilities for this $r\in\{1,2,\ldots,\ell\}\backslash\{a+m+1\}$.
We see, using
(\ref{gammaregeq}), (\ref{zetaregeq}), and (\ref{zetasingeq}),
that the resulting contribution to $\mathscr{C}_{a,0}^{[\ell]}$ is
\[
\sum_{\substack{j=1 \\ j\neq a+m+1}}^{\ell}
\zeta'(-a-m+j)\times \prod_{\substack{r=1 \\ r\notin\{a+m+1,j\}}}^{\ell}
\zeta(-a-m+r)\ .
\]
Adding all these contributions, while noting the cancellation due to $\Gamma'(1)=-\gamma$, and cleaning up, gives the formula stated in the proposition.
\qed

The last proposition will be further specialized, in the next two lemmas, to the most important $Z$ functions in this article, i.e., $Z_{\ell}^{[\ell]}$ and $Z_{\ell-1}^{[\ell]}$.
For the moment, we will leave the zeta functions unevaluated, so as to make the patterns present visually obvious.

\begin{lemma}
The coefficients of the complete asymptotic expansion of $Z_{\ell}^{[\ell]}(t)$
are given by
\[
\mathscr{C}_{a,b}^{[\ell]}(\ell,\ldots,1)=
\left\{
\begin{array}{ll}
(\ell-1)!\zeta(2)\zeta(3)\cdots\zeta(\ell)\ , & {\rm if}\ (a,b)=(-\ell,0)\ , \\
\zeta(0)(\ell-2)!\zeta(2)\zeta(3)\cdots\zeta(\ell-1)\ , & {\rm if}\ (a,b)=(-\ell+1,0)\ , \\
\zeta(-1)\zeta(0)(\ell-3)!\zeta(2)\zeta(3)\cdots\zeta(\ell-2)\ , & {\rm if}\ (a,b)=(-\ell+2,0)\ , \\
\vdots & \vdots \\
\zeta(-\ell+2)\zeta(-\ell+3)\cdots\zeta(-1)\zeta(0) 0! & {\rm if}\ (a,b)=(-1,0)\ , \\
\frac{(-1)^a}{a!}\zeta(-a-\ell+1)\zeta(-a-\ell+2)\cdots\zeta(-a)
\ , & {\rm if}\ a\ge 0\ {\rm and}\ b=0\ , \\
0\ & {\rm else}\ .
\end{array}
\right.
\]
\end{lemma}

\begin{lemma}
The coefficients of the complete asymptotic expansion of $Z_{\ell-1}^{[\ell]}(t)$
are given by
\begin{eqnarray*}
\lefteqn{
\mathscr{C}_{a,b}^{[\ell]}(\ell-1,\ldots,0)=} & & \\
 & & \left\{
\begin{array}{ll}
(\ell-2)!\zeta(2)\zeta(3)\cdots\zeta(\ell-1)\zeta(\ell)\ , & {\rm if}\ (a,b)=(-\ell+1,0)\ , \\
\zeta(0)(\ell-3)!\zeta(2)\zeta(3)\cdots\zeta(\ell-2)\zeta(\ell-1)\ , & {\rm if}\ (a,b)=(-\ell+2,0)\  , \\
\zeta(-1)\zeta(0)(\ell-4)!\zeta(2)\zeta(3)\cdots\zeta(\ell-3)\zeta(\ell-2)\ , & {\rm if}\ (a,b)=(-\ell+3,0)\ , \\
\vdots & \vdots \\
\zeta(-\ell+3)\zeta(-\ell+4)\cdots\zeta(-1)\zeta(0) 0! \zeta(2)\ , & {\rm if}\ (a,b)=(-1,0)\ , \\
\zeta(-\ell+2)\zeta(-\ell+3)\cdots\zeta(-2)\zeta(-1)\zeta(0)\ ,  & {\rm if}\ (a,b)=(0,1)\ , \\
\left.\frac{{\rm d}}{{\rm d}x}\zeta(x-\ell+2)\zeta(x-\ell+3)\cdots\zeta(x-2)\zeta(x-1)\zeta(x)\right|_{x=0}\ , & {\rm if}\ (a,b)=(0,0)\ , \\
\frac{(-1)^{a}}{a!}\zeta(-a-\ell+2)\zeta(-a-\ell+3)\cdots \zeta(-a+1)
\ , & {\rm if}\ a\ge 1\ {\rm and}\ b=0\ , \\
0\ & {\rm else}\ .
\end{array}
\right.
\end{eqnarray*}
\end{lemma}
No proof is needed since these two lemmas follow immediately from the last proposition. 
Now taking advantage of the trivial zeros of the zeta function $\zeta(-2)=\zeta(-4)=\zeta(-6)=\cdots=0$, we can further simplify the above asymptotic expansions as follows. Note that Bringmann, Franke and Heim also obtained such asymptotics, see~\cite[Eq. 3.5 and 3.6]{BringmannFH}.
For $\ell\ge 3$, we have
\begin{eqnarray}
Z_{\ell}^{[\ell]}(t)&=& (\ell-1)!\zeta(2)\zeta(3)\cdots\zeta(\ell)\ t^{-\ell}\nonumber \\
 & & +\zeta(0)(\ell-2)!\zeta(2)\zeta(3)\cdots\zeta(\ell-1)\ t^{-\ell+1}\nonumber \\
 & & +\zeta(-1)\zeta(0)(\ell-3)!\zeta(2)\cdots \zeta(\ell-2)\ t^{-\ell+2}\nonumber \\
 & & +O(t^{\infty})\ ,
\label{Zellasymeq}
\end{eqnarray}
while for $\ell=2$, we have
\begin{equation}
Z_{2}^{[2]}(t)=\zeta(2)t^{-2}+\zeta(0)t^{-1}+\zeta(-1)\zeta(0)
+O(t^{\infty})\ .
\label{Z22eq}
\end{equation}
For $\ell\ge 4$, we have
\begin{eqnarray}
Z_{\ell-1}^{[\ell]}(t)&=& (\ell-2)!\zeta(2)\zeta(3)\cdots\zeta(\ell)\ t^{-\ell+1}\nonumber \\
 & & +\zeta(0)(\ell-3)!\zeta(2)\zeta(3)\cdots\zeta(\ell-1)\ t^{-\ell+2}\nonumber \\
 & & +\zeta(-1)\zeta(0)(\ell-4)!\zeta(2)\cdots \zeta(\ell-2)\ t^{-\ell+3}\nonumber \\
 & & +\mathscr{C}_{0,0}^{[\ell]}(\ell-1,\ell-2,\ldots,0)\nonumber \\
 & & +O(t^{\infty})\ ,
\label{Zminusasymeq}
\end{eqnarray}
where the constant term reduces to
\[
\mathscr{C}_{0,0}^{[\ell]}(\ell-1,\ell-2,\ldots,0)=
\left\{
\begin{array}{ccc}
\zeta'(-2)\zeta(-1)\zeta(0) & {\rm if} & \ell=4\ , \\
\zeta(-3)\zeta'(-2)\zeta(-1)\zeta(0) & {\rm if} & \ell=5\ . \\
0 & {\rm if} & \ell\ge 6\ \ .
\end{array}
\right.
\]
For $\ell=2$, we have
\begin{eqnarray}
Z_{1}^{[2]}(t) & = & \zeta(2) t^{-1}+\zeta(0)(-\ln t)+\zeta'(0)-\zeta(-1)\zeta(0)t+O(t^{\infty}) \label{Z21eq} \\
 & = & \frac{\pi^2}{6}t^{-1}+\frac{1}{2}\ln t-\frac{1}{2}\ln(2\pi)-\frac{1}{24}t+O(t^{\infty}) \ , \nonumber
\end{eqnarray}
after evaluating the zeta functions. We thus recover a well known result (see~\cite[Exercise 3, p. 58]{deBruijn} or~\cite[Exercise 9, p. 10]{TenenbaumBook}).
Finally, for $\ell=3$, we have
\begin{eqnarray}
Z_{2}^{[3]}(t) & = & \zeta(2)\zeta(3) t^{-2}+\zeta(0) \zeta(2)  t^{-1}+ \zeta(-1) \zeta(0) (-\ln t) \nonumber \\
 & & +\left[\zeta'(-1) \zeta(0) +\zeta(-1) \zeta'(0) \right]+O(t^{\infty})
\label{Z32eq}
\end{eqnarray}
which we will not further evaluate, given the presence of the unwieldy $\zeta(3)$
and $\zeta'(-1)$.

The following functions will be needed later
\begin{equation}
U(t):=\mathcal{K}_{\ell}^{-1}t^{\ell}Z_{\ell}^{[\ell]}(t)\ ,
\label{Udefeq}
\end{equation}
and
\begin{equation}
V(t):=(\ell-1)\mathcal{K}_{\ell}^{-1}t^{\ell-1}Z_{\ell-1}^{[\ell]}(t)\ .
\label{Vdefeq}
\end{equation}
From (\ref{Zellasymeq}) and (\ref{Zminusasymeq}) we see that when $\ell\ge 4$, these admit the complete asymptotic expansions
\[
U(t)\sim A_0+A_1 t+A_2 t^2+A_3 t^3+\cdots
\]
and
\[
V(t)\sim B_0+B_1 t+B_2 t^2+B_3 t^3+\cdots
\]
where the $A$ and $B$ constants are the ones which have been defined in \S\ref{constdefsec}.
We note the following asymptotics, as $t\rightarrow 0^{+}$, which were given in~\cite[\S2]{AbdesselamS}, but also follow from Propostion \ref{Zmprop}:
\begin{eqnarray}
Z_{\ell+1}^{[\ell]}(t) & \sim & \ell\mathcal{K}_{\ell} t^{-\ell-1}\ , \label{Zplus1asymeq} \\
Z_{\ell+2}^{[\ell]}(t) & \sim & \ell(\ell+1)\mathcal{K}_{\ell} t^{-\ell-2}\ . \label{Zplus2asymeq} 
\end{eqnarray}
The next important property we will establish for the functions $Z_{m}^{[\ell]}(t)$ is their strict log-convexity, for fixed $t>0$, with respect to $m\in\mathbb{Z}$. This is the object of the next lemma.

\begin{lemma}\label{convexlem}
For all $t>0$ and all  $m\in\mathbb{Z}$, we have
\[
Z_{m}^{[\ell]}(t)^2<Z_{m-1}^{[\ell]}(t)\ Z_{m+1}^{[\ell]}(t)\ .
\]
\end{lemma}

\noindent{\bf Proof:}
By introducing a second set of summation indices $\eta_1,\ldots,\eta_{\ell}$ and expanding the products we get
\begin{eqnarray*}
\lefteqn{Z_{m+1}^{[\ell]}(t)Z_{m-1}^{[\ell]}(t)-Z_{m}^{[\ell]}(t)^2 =} & & \\
& & \sum_{\substack{\delta_1,\ldots,\delta_{\ell}\ge 1 \\ \eta_1,\ldots,\eta_{\ell}\ge 1}}
\left[\ \delta_1^{m}\delta_2^{m-1}\cdots\delta_{\ell}^{1}\ e^{-t\delta_1\cdots\delta_{\ell}}
\ \eta_1^{m-2}\eta_2^{m-3}\cdots\eta_{\ell}^{-1}\ e^{-t\eta_1\cdots\eta_{\ell}} 
\right. \\
 & & \qquad\qquad 
-\ \delta_1^{m-1}\delta_2^{m-2}\cdots\delta_{\ell}^{0}\ e^{-t\delta_1\cdots\delta_{\ell}}
\ \eta_1^{m-1}\eta_2^{m-2}\cdots\eta_{\ell}^{0}\ e^{-t\eta_1\cdots\eta_{\ell}}
\left.
\ \right] \\
 &= &  \sum_{\substack{\delta_1,\ldots,\delta_{\ell}\ge 1 \\ \eta_1,\ldots,\eta_{\ell}\ge 1}}
\delta_1^{m-2}\delta_2^{m-3}\cdots\delta_{\ell}^{-1}
\ \eta_1^{m-2}\eta_2^{m-3}\cdots\eta_{\ell}^{-1}
\ e^{-t(\delta_1\cdots\delta_{\ell}+\eta_1\cdots\eta_{\ell})}
\left[\delta_1^2\cdots\delta_{\ell}^2-\delta_1\cdots\delta_{\ell}
\eta_1\cdots\eta_{\ell}
\right] \\
 & = & \frac{1}{2}
\sum_{\substack{\delta_1,\ldots,\delta_{\ell}\ge 1 \\ \eta_1,\ldots,\eta_{\ell}\ge 1}}
\delta_1^{m-2}\delta_2^{m-3}\cdots\delta_{\ell}^{-1}
\ \eta_1^{m-2}\eta_2^{m-3}\cdots\eta_{\ell}^{-1}
\ e^{-t(\delta_1\cdots\delta_{\ell}+\eta_1\cdots\eta_{\ell})}
\left(\delta_1\cdots\delta_{\ell}-\eta_1\cdots\eta_{\ell}\right)^2 \\
 & > & 0\ .
\end{eqnarray*}
The last equality is obtained from the previous one by symmetrizing over the dummy tuples of summation indices $(\delta_1,\ldots,\delta_{\ell})$ and $(\eta_1,\ldots,\eta_{\ell})$.
The strict inequality follows by picking say all the $\delta$'s to be $2$ and all the $\eta$'s to be $1$.
\qed

\begin{remark}
It would be interesting to determine the sign of
\[
Z_{\alpha}^{[\ell]}(t)\ Z_{\alpha+e_i+e_j}^{[\ell]}(t)
-Z_{\alpha+e_i}^{[\ell]}(t)\ Z_{\alpha+e_j}^{[\ell]}(t)
\]
when $i\neq j$.
Here we used the vector notation $\alpha=(\alpha_1,\ldots,\alpha_{\ell})\in\mathbb{R}^{\ell}$  , and the $e_i$'s are the canonical basis vectors of $\mathbb{R}^{\ell}$. See~\cite{AbdesselamCorr} for a study of similar problems where the above symmetrization trick plays an important role.
\end{remark}

Consider now the $C^{\infty}$ function $h_{\ell}(t)$ defined in (\ref{hdefeq}).
By Lemma \ref{convexlem}, we have
\[
h'_{\ell}(t)=(\ell-1)\times\frac{Z_{\ell+1}^{[\ell]}(t)Z_{\ell-1}^{[\ell]}(t)-Z_{\ell}^{[\ell]}(t)^2}{Z_{\ell}^{[\ell]}(t)^2}>0\ ,
\]
and the function $h_{\ell}$ is strictly increasing.
By the asymptotics (\ref{Zellasymeq}) and (\ref{Zminusasymeq}), we have $h_{\ell}(t)\sim t$ when $t\rightarrow 0^{+}$ and thus
$\lim_{t\rightarrow 0^{+}}h_{\ell}(t)=0$.
We leave it as an easy exercise for the reader to check that for any $m\in\mathbb{Z}$, we have $Z_{m}^{[\ell]}(t)\sim e^{-t}$ when $t\rightarrow\infty$. This implies $\lim_{t\rightarrow\infty}h_{\ell}(t)=\ell-1$. Putting all this together, we see that $h_{\ell}$ is an increasing $C^{\infty}$ diffeomorphism from $(0,\infty)$ onto $(0,\ell-1)$.
The choices made in (\ref{tndefeq}) and (\ref{rhondefeq}) are now well defined.
Since $h_{\ell}(t)\sim t$ and $k_n\sim sn^{\frac{\ell-1}{\ell}}$ by hypothesis, it is easy to see that, as $n\rightarrow\infty$, we have
\begin{equation}
t_n \sim (\ell-1) s n^{-\frac{1}{\ell}}\ .
\label{tnasymeq}
\end{equation}
From (\ref{tnasymeq}) and from the leading asymptotic in (\ref{Zellasymeq}), we deduce
\begin{equation}
\lim\limits_{n\rightarrow\infty} \rho_n=\frac{(\ell-1)^{\ell} s^{\ell}}{\mathcal{K}_{\ell}}=:\rho_{\infty}\ .
\label{rhonlimeq}
\end{equation}
In preparation for the analysis in \S\ref{saddlesec}, we also define the quantities
\begin{eqnarray}
\lambda_n & := & Z_{\ell+1}^{[\ell]}(t_n)^{-\frac{1}{2}}\ , \label{lambdadefeq} \\
\mu_n & := & Z_{\ell-1}^{[\ell]}(t_n)^{-\frac{1}{2}}\ , \label{mudefeq}
\end{eqnarray}
which will later be used to do a change of variables $\theta=\lambda_n\Theta$, and $\omega=\mu_n\Omega$.
By (\ref{Zplus1asymeq}) and (\ref{Zminusasymeq}), these quantities obey the following asymptotic equivalences
\begin{eqnarray}
\lambda_n & \sim & \frac{1}{\sqrt{\ell\mathcal{K}_{\ell}}}\times t_n^{\frac{\ell+1}{2}}\ , \label{lambdaasymeq} \\
\mu_n & \sim & \sqrt{\frac{\ell-1}{\mathcal{K}_{\ell}}}\times t_n^{\frac{\ell-1}{2}}\ .
\label{muasymeq}
\end{eqnarray}

\section{Bivariate saddle point analysis}
\label{saddlesec}

Picking up the thread from \S\ref{outlinesec}, we can write the integrand in (\ref{Indefeq}) as
\[
\mathcal{I}_{n}(\theta,\omega)=\exp(-Q_n(\theta,\omega))\ ,
\]
with
\begin{equation}
Q_n(\theta,\omega):=
in\theta+i k_n\omega+\rho_n
\sum_{\delta_1,\ldots,\delta_{\ell=1}}^{\infty}
\delta_{1}^{\ell-2}\cdots\delta_{\ell}^{-1}
\ e^{-\delta_1\cdots\delta_{\ell} t_n}\left(
1-e^{i(\omega+\delta_1\cdots\delta_{\ell}\theta)}
\right)\ .
\label{Qdefeq}
\end{equation}
In this section, we will not separate the cases where $\ell\ge 4$ and $\ell\in\{2,3\}$.
We will break the domain of integration into two regions: the major arc or two-dimensional patch rather,
\[
\mathcal{D}_{\rm maj}:=\left\{
(\theta,\omega)\in(-\pi,\pi)^2\ |\ |\theta|\le t_n\ {\rm and}\ |\omega|\le\frac{\pi}{2}
\right\}\ ,
\]
and the minor patch
\[
\mathcal{D}_{\rm min}:=(-\pi,\pi)^2\backslash
\mathcal{D}_{\rm maj}\ .
\]
Proceeding similarly to~\cite[\S3]{AbdesselamS}, we note that
\begin{eqnarray}
{\rm Re}\ Q_n(\theta,\omega) & =&
2\rho_n\sum_{\delta_1,\ldots,\delta_{\ell=1}}^{\infty}
\delta_{1}^{\ell-2}\cdots\delta_{\ell}^{-1}
\ e^{-\delta_1\cdots\delta_{\ell} t_n}
\sin^2\left(\frac{\omega+\delta_1\cdots\delta_{\ell}\theta}{2}\right) \nonumber \\
 & \ge & 2\rho_n\sum_{p=1}^{\infty}
p^{\ell-2}\ e^{-p t_n} \sin^2\left(\frac{\omega+p\theta}{2}\right)\ ,
\label{Qlowbdeq}
\end{eqnarray}
where we exploited positivity, and only kept tuples of the form $(\delta_1,\ldots,\delta_{\ell})=(p,1,\ldots,1)$.

\subsection{Major patch estimates}

In this section we will provide estimates for when $(\theta,\omega)\in\mathcal{D}_{\rm maj}$.
We define the cutoff
\[
N_n:=\left\lfloor\frac{\pi}{2 t_n}\right\rfloor\ ,
\]
and use it to weaken the lower bound (\ref{Qlowbdeq}) to
\[
{\rm Re}\ Q_n(\theta,\omega) \ge
2\rho_n\sum_{p=1}^{N_n}
p^{\ell-2}\ e^{-p t_n} \sin^2\left|\frac{\omega+p\theta}{2}\right|\ .
\]
For $1\le p\le N_n$, and $|\theta|\le t_n$, we have $|p\theta|\le\frac{\pi}{2}$. We also have $|\omega|\le\frac{\pi}{2}$ by hypothesis. Hence $\left|\frac{\omega+p\theta}{2}\right|\le\frac{\pi}{2}$.
We now use the concavity of the sine function on the interval $[0,\frac{\pi}{2}]$, i.e., the bound $\sin u\ge \frac{2u}{\pi}$ for $u$ in that interval. 
This gives
\begin{eqnarray*}
{\rm Re}\ Q_n(\theta,\omega) & \ge &
2\rho_n\sum_{p=1}^{N_n}
p^{\ell-2} e^{-p t_n} \left(\frac{2}{\pi}
\left|\frac{\omega+p\theta}{2}\right|\right)^2 \\
 & \ge & \frac{2\rho_n}{\pi^2 e^{\frac{\pi}{2}}}
\sum_{p=1}^{N_n}p^{\ell-2}(\omega+p\theta)^2\ ,
\end{eqnarray*}
where we also used $pt_n\le N_n t_n\le\frac{\pi}{2}$.
By expanding the square, we can rewrite the previous inequality as
\[
{\rm Re}\ Q_n(\theta,\omega)  \ge 
\frac{2\rho_n}{\pi^2 e^{\frac{\pi}{2}}}
\left[
u_n(\lambda_{n}^{-1}\theta)^2+2v_n(\lambda_{n}^{-1}\theta)(\mu_n^{-1}\omega)+w_n(\mu_{n}^{-1}\omega)^2
\right]\ ,
\]
where we introduced the new quantities
\begin{eqnarray*}
u_n & := & \lambda_n^2\times\sum_{p=1}^{N_n} p^{\ell} \ ,\\
v_n & := & \lambda_n\mu_n\times\sum_{p=1}^{N_n} p^{\ell-1} \ ,\\
w_n & := & \mu_n^2\times\sum_{p=1}^{N_n} p^{\ell-2} \ .\\
\end{eqnarray*}
For any $m\ge 0$, the elementary inequalities
\[
\int_0^N x^m\ {\rm d}x\le\sum_{p=1}^{N}p^m\le
\int_1^{N+1} x^m\ {\rm d}x
\]
immediately give the asymptotic equivalence
\[
\sum_{p=1}^{N}p^m\sim\frac{N^{m+1}}{m+1}\ ,
\]
when $N\rightarrow\infty$, which we will use for $m=\ell,\ell-1,\ell-2$, together with the trivial $N_n\sim\frac{\pi}{2t_n}$.
This combined with (\ref{lambdaasymeq}) and (\ref{muasymeq}) results in the computation of limits:
\begin{eqnarray*}
\lim\limits_{n\rightarrow\infty}u_n & = & \frac{1}{\ell(\ell+1)\mathcal{K}_{\ell}}
\times\left(\frac{\pi}{2}\right)^{\ell+1} =: u_{\infty}\ , \\
\lim\limits_{n\rightarrow\infty}v_n & = & \frac{1}{\mathcal{K}_{\ell}}
\times\sqrt{\frac{\ell-1}{\ell^3}}\times
\left(\frac{\pi}{2}\right)^{\ell} =: v_{\infty} \ , \\
\lim\limits_{n\rightarrow\infty}w_n & = & \frac{1}{\mathcal{K}_{\ell}}
\times\left(\frac{\pi}{2}\right)^{\ell-1} =: w_{\infty}\ . \\
\end{eqnarray*}
We note that $u_{\infty}>0$, $w_{\infty}>0$, and also, with some easy algebra,
\[
u_{\infty}w_{\infty}-v_{\infty}^2=\frac{1}{\ell^3(\ell+1)\mathcal{K}_{\ell}^2}
\times\left(\frac{\pi}{2}\right)^{2\ell}
>0\ .
\]
This means that the matrix
\[
\begin{pmatrix}
u_{\infty} & v_{\infty} \\
v_{\infty} & w_{\infty}
\end{pmatrix}
\]
belongs to the cone ${\rm Pos}_2$ 
of positive definite real symmetric $2\times 2$ matrices. 
Since ${\rm Pos}_2$ is open in the space of real symmetric $2\times 2$ matrices, one can pick $\epsilon>0$ such that
\[
\begin{pmatrix}
u_{\infty}-\epsilon & v_{\infty} \\
v_{\infty} & w_{\infty}-\epsilon
\end{pmatrix}
\in{\rm Pos}_2\ .
\]
Again, since ${\rm Pos}_2$ is open, and because of the above limits, we have that
\[
\begin{pmatrix}
u_{n}-\epsilon & v_{n} \\
v_{n} & w_{n}-\epsilon
\end{pmatrix}
\in{\rm Pos}_2\ ,
\]
for $n$ large enough. For such $n$, we have $\forall (\Theta,\Omega)\in\mathbb{R}^2$,
\[
u_n\Theta^2+2v_n\Theta\Omega+w_n\Omega^2\ge
\epsilon(\Theta^2+\Omega^2)\ . 
\]
By applying this to $\Theta=\lambda_{n}^{-1}\theta$ and $\Omega=\mu_{n}^{-1}\omega$, we get the lower bound
\[
{\rm Re}\ Q_n(\theta,\omega)\ge
\frac{2\rho_n}{\pi^2 e^{\frac{\pi}{2}}}\ \epsilon
\left(\lambda_{n}^{-2}\theta^2+\mu_{n}^{-2}\omega^2\right)\ .
\]
Now pick some $\eta_{\rm maj}>0$ such that $\eta_{\rm maj}<
\frac{2\rho_{\infty}}{\pi^2 e^{\frac{\pi}{2}}}\epsilon$. We have established that for $n$ large enough, $\forall(\theta,\omega)\in\mathcal{D}_{\rm maj}$ we have
\begin{equation}
{\rm Re}\ Q_n(\theta,\omega)\ge
\eta_{\rm maj}
\left(\lambda_{n}^{-2}\theta^2+\mu_{n}^{-2}\omega^2\right)\ .
\label{Qmajbdeq}
\end{equation}

\subsection{Minor patch estimates}

In this section we will provide estimates for when $(\theta,\omega)\in\mathcal{D}_{\rm min}$. We go back to (\ref{Qlowbdeq}), use the inequality $p^{\ell-2}\ge 1$ and reverse the use of the identity ${\rm Re}(1-e^{iu})=2\sin^2\left(\frac{u}{2}\right)$.
We now have the lower bound
\[
{\rm Re}\ Q_n(\theta,\omega)\ge 2\rho_n
\ {\rm Re}\ W_n(\theta,\omega)\ ,
\]
where
\begin{eqnarray*}
W_n(\theta,\omega) & := &
\sum_{p=1}^{\infty} e^{-pt_n}
\left(1-e^{i\omega+ip\theta}\right) \\
 & = & \frac{e^{-t_n}}{1-e^{-t_n}}-e^{i\omega}\times\frac{e^{-t_n+i\theta}}{1-e^{-t_n+i\theta}} \\
 & = & \frac{1}{e^{t_n}-1}-\frac{e^{i\omega}}{e^{t_n-i\theta}-1}\ .
\end{eqnarray*}
Taking the real part, we see that
\begin{eqnarray*}
{\rm Re}\ W_n(\theta,\omega) & = &
\frac{1}{e^{t_n}-1}-{\rm Re}\left[
\frac{e^{i\omega}\left(e^{t_n+i\theta}-1\right)}{\left(e^{t_n-i\theta}-1\right)
\left(e^{t_n+i\theta}-1\right)}
\right] \\
 & = & \frac{1}{e^{t_n}-1}+\frac{X_n(\theta,\omega)}{Y_n(\theta)}\ ,
\end{eqnarray*}
where
\[
X_n(\theta,\omega):=\cos\omega-e^{t_n}\cos(\omega+\theta)\ ,
\]
and
\[
Y_n(\theta):=e^{2t_n}-2e^{t_n}\cos\theta +1>0\ .
\]
By expanding the cosine, and comparing $\cos \theta$ to $1$, we can write
\begin{eqnarray*}
X_n(\theta,\omega) & = & \cos\omega-e^{t_n}\cos\omega\cos\theta
+e^{t_n}\sin\omega\sin\theta \\
 & = & (-\cos\omega)\left(e^{t_n}-1\right)+
e^{t_n}\cos\omega(1-\cos\theta)+e^{t_n}\sin\omega\sin\theta\ .
\end{eqnarray*}
We note that
\[
\left|\frac{e^{t_n}\cos\omega(1-\cos\theta)}{Y_n(\theta)}\right|\le 
\frac{e^{t_n}(1-\cos\theta)}{Y_n(\theta)}
=\frac{e^{t_n}(1-\cos\theta)}{(e^{t_n}-1)^2+2e^{t_n}(1-\cos\theta)}\ .
\]
By dropping the term with the square, in the denominator, we immediately obtain
\[
\left|\frac{e^{t_n}\cos\omega(1-\cos\theta)}{Y_n(\theta)}\right|\le\frac{1}{2}\ ,
\]
which results in the key inequality
\begin{equation}
{\rm Re}\ W_n(\theta,\omega)\ge
\frac{1}{e^{t_n}-1}-\frac{1}{2}+\frac{1}{Y_n(\theta)}
\left[(-\cos\omega)(e^{t_n}-1)+e^{t_n}\sin\omega\sin\theta\right]\ .
\label{keymineq}
\end{equation}

\noindent\underline{1st case:} We assume $\omega\in(-\pi,\pi)$ satisfies $|\omega|>\frac{\pi}{2}$. Then, $-\cos\omega\ge 0$, and we simply toss the corresponding term in (\ref{keymineq}). We use $|\sin\omega|\le 1$ in the last term, and thus get the lower bound
\[
{\rm Re}\ W_n(\theta,\omega)\ge
\frac{1}{e^{t_n}-1}-\frac{1}{2}-\frac{e^{t_n}\sin\theta}{Y_n(\theta)}\ .
\]
Now consider the calculus exercise of determining the maximum
\[
\max_{-\pi\le\theta\le\pi}|f(\theta)|=\max_{0\le \theta\le\pi} f(\theta)
\]
for the odd function
\[
f(\theta)=\frac{\sin\theta}{c-\cos\theta}\ ,
\]
where $c>1$ is some constant.
The derivative is
\[
f'(\theta)=\frac{\cos\theta -1}{(c-\cos\theta)^2}\ .
\]
Therefore, the maximum is reached at $\theta=\arccos(c^{-1})$, and we have
\[
\max_{-\pi\le\theta\le\pi}|f(\theta)|=\frac{\sqrt{1-\left(\frac{1}{c}\right)^2}}{c-\frac{1}{c}}=\frac{1}{\sqrt{c^2-1}}\ .
\]
We use this for $c=\frac{e^{2t_n}+1}{2e^{t_n}}>1$, and deduce
\[
\left|
\frac{\sin\theta}{Y_n(\theta)}
\right|\le\frac{1}{2e^{t_n}}\times\frac{1}{\sqrt{\left(
\frac{e^{2t_n}+1}{2e^{t_n}}
\right)^2-1}}=
\frac{1}{e^{2t_n}-1}\ .
\]
Hence, we obtain
\begin{equation}
{\rm Re}\ W_n(\theta,\omega)\ge
\frac{1}{e^{t_n}-1}-\frac{1}{2}-\frac{e^{t_n}}{e^{2t_n}-1}
=\frac{1}{e^{2t_n}-1}-\frac{1}{2}\sim\frac{1}{2t_n}\ .
\label{mincase1eq}
\end{equation}

\noindent\underline{2nd case:} We assume $\theta\in(-\pi,\pi)$ satisfies $|\theta|>t_n$.
For any real numbers $c_1,c_2$ one can of course write, for some suitable angle $\tau$,
\[
c_1\cos\omega+c_2\sin\omega=\sqrt{c_1^2+c_2^2}\times\cos(\omega-\tau)\ge-
\sqrt{c_1^2+c_2^2}\ .
\]
By applying that to $c_1=1-e^{t_n}$ and $c_2=e^{t_n}\sin\theta$,
we deduce from (\ref{keymineq}) that
\[
{\rm Re}\ W_n(\theta,\omega)\ge
\frac{1}{e^{t_n}-1}-\frac{1}{2}
-\frac{1}{Y_n(\theta)}\sqrt{
(e^{t_n}-1)^2+e^{2t_n}\sin^2\theta}\ .
\]
We now notice that
\begin{eqnarray*}
\sqrt{
(e^{t_n}-1)^2+e^{2t_n}\sin^2\theta} & = &
e^{\frac{t_n}{2}}\sqrt{
e^{-t_n}(e^{t_n}-1)^2+e^{t_n}(1-\cos\theta)(1+\cos\theta)} \\
 & \le & e^{\frac{t_n}{2}}\sqrt{
(e^{t_n}-1)^2+2e^{t_n}(1-\cos\theta)} \ ,
\end{eqnarray*}
with the expression in the last square root being an alternate formula for $Y_n(\theta)$.
As a result,
\[
{\rm Re}\ W_n(\theta,\omega)\ge
\frac{1}{e^{t_n}-1}-\frac{1}{2}
-\frac{e^{\frac{t_n}{2}}}{\sqrt{Y_n(\theta)}}\ .
\]
Now for
$t_n<|\theta|<\pi$ we have $\cos\theta\le\cos t_n$ and therefore
\[
Y_n(\theta)\ge e^{2t_n}-2 e^{t_n}\cos t_n+1=2t_n^2+O(t_n^3)\ ,
\]
by an asymptotic computation done in~\cite[\S4.2]{AbdesselamS}.
We have thus shown, in this second case, that 
\begin{equation}
{\rm Re}\ W_n(\theta,\omega)\ge
\frac{1}{e^{t_n}-1}-\frac{1}{2}
-\frac{e^{\frac{t_n}{2}}}{\sqrt{e^{2t_n}-2 e^{t_n}\cos t_n+1}}
\sim\frac{1-\frac{1}{\sqrt{2}}}{t_n}\ .
\label{mincase2eq}
\end{equation}

We now combine both cases and the corresponding lower bounds (\ref{mincase1eq}) and (\ref{mincase2eq}), by picking some $\epsilon>0$, such that $\epsilon<1-\frac{1}{\sqrt{2}}<\frac{1}{2}$. For $n$ large enough we have 
$\forall (\theta,\omega)\in\mathcal{D}_{\rm min}$,
\[
{\rm Re}\ W_n(\theta,\omega)\ge \epsilon\ t_n^{-1}\ .
\]
Finally, pick some $\eta_{\rm min}>0$ such that $\eta_{\rm min}<\epsilon\times\frac{[(\ell-1)s]^{\ell-1}}{\mathcal{K}_{\ell}}$.
From the value of $\rho_{\infty}$ in (\ref{rhonlimeq}) and the asymptotics (\ref{tnasymeq}), we easily see that, 
for $n$ large enough, we have 
$\forall (\theta,\omega)\in\mathcal{D}_{\rm min}$,
\[
{\rm Re}\ Q_n(\theta,\omega)\ge \eta_{\rm min}\ n^{\frac{1}{\ell}}\ .
\]
This gives us, for $n$ large enough, the upper bound
\begin{equation}
\left|\int_{\mathcal{D}_{\rm min}} e^{-Q_n(\theta,\omega)}
\ \frac{{\rm d}\theta{\rm d}\omega}{(2\pi)^2}
\right|\le \exp\left(-\eta_{\rm min}\ n^{\frac{1}{\ell}}\right)\ .
\label{minintbdeq}
\end{equation}

\subsection{Proof of the asymptotics for the $\mathscr{I}$ integral}

We now perform the change of variables $\theta=\lambda_n\Theta$, and $\omega=\mu_n\Omega$, and write the main integral as
\[
\int_{\mathcal{D}_{\rm maj}} e^{-Q_n(\theta,\omega)}
\ \frac{{\rm d}\theta{\rm d}\omega}{(2\pi)^2}
=\frac{\lambda_n\mu_n}{(2\pi)^2}
\int_{\mathbb{R}^2} f_n(\Theta,\Omega)\ {\rm d}\Theta{\rm d}\Omega\ ,
\]
where
\[
f_n(\Theta,\Omega):=\bbone\left\{
|\Theta|\le t_n\lambda_{n}^{-1}\ {\rm and}\ |\Omega|\le\frac{\pi}{2}\mu_{n}^{-1}\right\}
\times
 e^{-Q_n(\lambda_n\Theta,\mu_n\Omega)}\ .
\]
We note that
\[
t_n\lambda_{n}^{-1}\sim
\sqrt{\ell\mathcal{K}_{\ell}}\times t_n^{-\left(\frac{\ell-1}{2}\right)}
\rightarrow\infty\ ,
\]
and
\[
\frac{\pi}{2}\mu_{n}^{-1}\sim
\frac{\pi}{2}\sqrt{\frac{\mathcal{K}_{\ell}}{\ell-1}}\times t_n^{-\left(\frac{\ell-1}{2}\right)}
\rightarrow\infty\ ,
\]
because $\ell\ge 2$, by hypothesis.
We also have, thanks to (\ref{Qmajbdeq}), 
$\forall(\Theta,\Omega)\in\mathbb{R}^2$,
\[
|f_n(\Theta,\Omega)|\le e^{-{\rm Re}\ Q_n(\lambda_n\Theta,\mu_n\Omega)}
\le\exp\left(-\eta_{\rm maj}(\Theta^2+\Omega^2)\right)
\]
which is integrable, and therefore puts us in a position to use the dominated convergence theorem.

As in~\cite[\S4.3]{AbdesselamS}, we let $R(w):=e^w-1-w-\frac{w^2}{2}$, which has the property that $|R(w)|\le\frac{|w|^3}{6}$ whenever $w$ is pure imaginary.
We use this for $w=i(\omega+\delta_1\cdots\delta_{\ell}\theta)$ inside formula (\ref{Qdefeq}). This produces the rewriting
\begin{eqnarray}
Q_n(\lambda_n\Theta,\mu_n\Omega) & = & 
in\lambda_n\Theta+ik_n\mu_n\Omega \nonumber \\
 & & -i\rho_n\lambda_n\Theta Z_{\ell}^{[\ell]}(t_n)
-i\rho_n\mu_n\Omega Z_{\ell-1}^{[\ell]}(t_n) \nonumber \\
 & & +\frac{\rho_n}{2}\sum_{\delta_1,\ldots,\delta_{\ell}=1}^{\infty}
\delta_{1}^{\ell-2}\cdots\delta_{\ell}^{-1}\ e^{-\delta_1\cdots\delta_{\ell}t_n}\ \left(\mu_n\Omega+\delta_1\cdots\delta_{\ell}\lambda_n\Theta\right)^2 \nonumber \\
 & & +\mathcal{E}\ , 
\label{bigQeq}
\end{eqnarray}
with the error term given by
\[
\mathcal{E}:=-\rho_n
\sum_{\delta_1,\ldots,\delta_{\ell}=1}^{\infty}
\delta_{1}^{\ell-2}\cdots\delta_{\ell}^{-1}\ e^{-\delta_1\cdots\delta_{\ell}t_n}
\ R(i\mu_n\Omega+i\delta_1\cdots\delta_{\ell}\lambda_n\Theta)\ .
\]
We have, by the previous bound on $R(w)$, and by expanding the cube
\begin{eqnarray*}
|\mathcal{E}| & \le & 
\frac{\rho_n}{6}\sum_{\delta_1,\ldots,\delta_{\ell}=1}^{\infty}
\delta_{1}^{\ell-2}\cdots\delta_{\ell}^{-1}\ e^{-\delta_1\cdots\delta_{\ell}t_n}\ \left(\mu_n|\Omega|+\delta_1\cdots\delta_{\ell}\lambda_n|\Theta|\right)^3 \\
 & \le & \frac{1}{6}\rho_n\mu_n^{3}|\Omega|^3\ Z_{\ell-1}^{[\ell]}(t_n) \\
 & & + \frac{1}{2}\rho_n\mu_n^{2}\lambda_n\Omega^2|\Theta|\ Z_{\ell}^{[\ell]}(t_n) \\
 & & +\frac{1}{2}\rho_n\mu_n\lambda_{n}^{2}|\Omega|\Theta^2\ Z_{\ell+1}^{[\ell]}(t_n) \\
 & &+\frac{1}{6}\rho_n\lambda_n^{3}|\Theta|^3\ Z_{\ell+2}^{[\ell]}(t_n) \ .
\end{eqnarray*}
From the limit (\ref{rhonlimeq}), and the asymptotics 
(\ref{Zminusasymeq}), (\ref{Zellasymeq}), (\ref{Zplus1asymeq}), (\ref{Zplus2asymeq}), (\ref{lambdaasymeq}), and (\ref{muasymeq}),
it is easy to see that all four terms in the last inequality have order $t_n^{\frac{\ell-1}{2}}\rightarrow 0$, because $\ell\ge 2$.
Hence, for fixed $\Theta,\Omega$, the error term $\mathcal{E}$ goes to zero, when $n\rightarrow\infty$.

The terms in (\ref{bigQeq}) which are linear in $\Theta,\Omega$ vanish by design, i.e., because of our choices (\ref{tndefeq}) and (\ref{rhondefeq}).
Once the square is expanded, the quadratic part in (\ref{bigQeq}) becomes
\[
\frac{1}{2}\rho_n\lambda_n^{2}\Theta^2\ Z_{\ell+1}^{[\ell]}(t_n)
+\rho_n\lambda_n\mu_n\Theta\Omega\ Z_{\ell}^{[\ell]}(t_n)  
+\frac{1}{2}\rho_n\mu_n^{2}\Omega^2\ Z_{\ell-1}^{[\ell]}(t_n) \ .
\]
From the limit (\ref{rhonlimeq}), and the asymptotics 
(\ref{Zminusasymeq}), (\ref{Zellasymeq}), (\ref{Zplus1asymeq}), (\ref{lambdaasymeq}), and (\ref{muasymeq}),
one easily finds that, for fixed $\Theta,\Omega$,
\begin{eqnarray*}
\lim\limits_{n\rightarrow\infty}
Q_n(\lambda_n\Theta,\mu_n\Omega) & = & \frac{(\ell-1)^{\ell}s^{\ell}}{\mathcal{K}_{\ell}}\times
\left[\frac{1}{2}\Theta^2+\sqrt{\frac{\ell-1}{\ell}}\ \Theta\Omega+\frac{1}{2}\Omega^2\right]\ \\
 & = & \frac{1}{2}(\Theta,\Omega) M (\Theta,\Omega)^{\rm T}\ ,
\end{eqnarray*}
where we introduced the matrix
\[
M:=\frac{(\ell-1)^{\ell}s^{\ell}}{\mathcal{K}_{\ell}}
\begin{pmatrix}
1 & \sqrt{\frac{\ell-1}{\ell}}\\
\sqrt{\frac{\ell-1}{\ell}} & 1
\end{pmatrix}\ .
\]
From the dominated convergence theorem, we finally get
\begin{eqnarray*}
\lim\limits_{n\rightarrow\infty}
\int_{\mathbb{R}^2} f_n(\Theta,\Omega)\ {\rm d}\Theta{\rm d}\Omega
 & = & \int_{\mathbb{R}^2} e^{-\frac{1}{2}(\Theta,\Omega) M (\Theta,\Omega)^{\rm T}}\ {\rm d}\Theta{\rm d}\Omega \\
 & = & \frac{2\pi}{\sqrt{{\rm det}\ M}} \\
 & = & \frac{2\pi\mathcal{K}_{\ell}\sqrt{\ell}}{(\ell-1)^{\ell}s^{\ell}}\ .
\end{eqnarray*}
By the asymptotics (\ref{lambdaasymeq}), (\ref{muasymeq}), and (\ref{tnasymeq}),
we have
\begin{eqnarray*}
\frac{\lambda_n\mu_n}{(2\pi)^2} & \sim & \frac{1}{(2\pi)^2}\times
\sqrt{\frac{\ell-1}{\ell}}\times\frac{1}{\mathcal{K}_{\ell}}\times t_{n}^{\ell} \\
 & \sim &  \frac{1}{(2\pi)^2}\times
\sqrt{\frac{\ell-1}{\ell}}\times\frac{1}{\mathcal{K}_{\ell}}\times
(\ell-1)^{\ell}s^{\ell}\times \frac{1}{n}\ .
\end{eqnarray*}
Putting things together we get
\[
\int_{\mathcal{D}_{\rm maj}} e^{-Q_n(\theta,\omega)}
\ \frac{{\rm d}\theta{\rm d}\omega}{(2\pi)^2}
\ \sim\ \frac{\sqrt{\ell-1}}{2\pi n}\ .
\]
Since this is a power law, it dominates the similar integral over
$\mathcal{D}_{\rm min}$
which satisfies the fractional exponential decay in (\ref{minintbdeq}).
As a result, we have established the asymptotic equivalence
\[
\mathscr{I}_n\sim \ \frac{\sqrt{\ell-1}}{2\pi n}\ ,
\]
for the integral defined in (\ref{curlyIdefeq}).

\section{Completion of the proof of the main theorem}

Since the analysis for the integral $\mathscr{I}_{n}$ in (\ref{multdeceq})
has been settled, we now turn our attention to the asymptotics of the prefactor $\mathscr{M}_{\ell}(t_n,\rho_n)$. We will produce asymptotics for the logarithm of this prefactor. Unlike the similar analysis done in~\cite{AbdesselamS} which only concerned the leading term, this time, we need the detailed
asymptotics of $\ln \mathscr{M}_{\ell}(t_n,\rho_n)$ all the way down to $o(1)$.
Using (\ref{logMeq}), (\ref{tndefeq}), and (\ref{rhondefeq}), we can write
\begin{equation}
\ln \mathscr{M}_{\ell}(t_n,\rho_n)=
nt_n-k_n\ln\rho_n
+k_n\ .
\label{betterlogMeq}
\end{equation}
Since we have
\begin{eqnarray*}
t_n & \sim & (\ell-1)\times\frac{k_n}{n}\ , \\
\rho_n & \sim & \frac{(\ell-1)^{\ell}}{\mathcal{K}_{\ell}}\times \frac{k_n^{\ell}}{n^{\ell-1}}\ ,
\end{eqnarray*}
this suggests
defining the small quantities $\varepsilon_n$, $\eta_n$ so that we have the equalities
\begin{eqnarray*}
t_n & = & (\ell-1)\times\frac{k_n}{n}\times(1+\varepsilon_n)\ , \\
\rho_n & = & \frac{(\ell-1)^{\ell}}{\mathcal{K}_{\ell}}\times \frac{k_n^{\ell}}{n^{\ell-1}}\times(1+\eta_n)\ .
\end{eqnarray*}
In other words, we let
\begin{eqnarray*}
\varepsilon_n & := & \frac{nt_n}{(\ell-1)k_n}-1\ , \\
\eta_n & := & \rho_n\times\frac{\mathcal{K}_{\ell}n^{\ell-1}}{(\ell-1)^{\ell}k_n^{\ell}}-1\ .
\end{eqnarray*}
We insert the formulas for $t_n$, $\rho_n$ in terms of $\varepsilon_n$, $\eta_n$, into (\ref{betterlogMeq}), and with a little bit of algebra, readily obtain
\[
\ln \mathscr{M}_{\ell}(t_n,\rho_n)=\mathscr{S}_{\ell}(n,k_n)+k_n\Delta_n\ ,
\]
where $\mathscr{S}_{\ell}$ is the expression discovered by Starr featuring in our main theorems,  and where
\[
\Delta_n:=(\ell-1)\varepsilon_n-\ln(1+\eta_n)\ .
\]
Before going into the derivation of the asymptotics for $\varepsilon_n$ and $\eta_n$, we pause in order to discuss the role of formal power series in asymptotic analysis. If $f(t)$ is some function defined in some domain $\mathcal{D}$ which has $0$ as a limit point, writing
\[
f(t)\sim a_0+a_1 t+a_2t^2+a_3 t^3+\cdots
\]
is defined, in every textbook on asymptotic analysis, as the statement that for all $N\ge 0$, the remainder $f(t)-(a_0+a_1t+a_2t^2+\cdots+a_N t^N)$ is $o(t^N)$ as $t$ goes to zero. What unfortunately is not often emphasized in textbooks, with the notable exception of the one by de Bruijn~\cite[\S1.6]{deBruijn},
is that the above statement is best understood as a precise relation $f(t)\sim A(t)$ between a true function on the left-hand side, and a formal power series
\[
A(t):=a_0+a_1t+a_2t^2+a_3 t^3+\cdots\ \in\ \mathbb{C}[[t]]\ ,
\]
on the right-hand side.
Each of two worlds connected by such a statement has its own natural internal operations. For functions, one defines sum, products, division, composition, etc.
in a pointwise manner for $t\in\mathcal{D}$. For formal power series one has natural analogues, since $\mathbb{C}[[t]]$ is a commutative ring, with invertible elements corresponding to formal power series with nonzero constant terms. One also has the operation of substitution (composition) of a formal power series without constant term into another formal power series. A proper handling of this operation mentioned in~\cite[\S1.6]{deBruijn} would need the notion of summable families and a topology on $\mathbb{C}[[t]]$ which can be seen as $\mathbb{C}^{\mathbb{N}}$.
The most used topology here is the product of discrete topologies on $\mathbb{C}$.
Perhaps the reason for the lack of emphasis on formal power series is that this is too much abstract algebra for analysis majors. Nevertheless, the beauty and computational effectiveness of the calculus of asymptotic power series
reflects the ``functorial'' nature of relations $f(t)\sim A(t)$ between functions and formal power series. For example, if $f(t)\sim A(t)$ and $g(t)\sim B(t)$, then $f(t)g(t)\sim A(t)B(t)$ where, by definition, the right-hand side is the product in the ring $\mathbb{C}[[t]]$.
The same holds for the other elementary operations, including composition~\cite[\S1.6]{deBruijn}. An exercise which we leave to the reader is that if $f(t)\sim t+a_2 t^2+a_3 t^3+\cdots$ then the locally well defined inverse function $f^{-1}$ also has an asymptotic expansion $t+b_2 t+b_3 t^3+\cdots$ where the last series is the compositional inverse of the formal power series $t+a_2 t^2+a_3 t^3+\cdots$. Determining this compositional inverse is the problem of reversion, and its explicit solution typically involves sums over tree graphs (see~\cite{AbdesselamJPA,AbdesselamSLC} and references therein).
Also note that formal power series like $t+a_2 t^2+a_3 t^3+\cdots$ form, for the operation of composition, a group one may call the group of unipotent formal diffeomorphisms.

We now put the above principles into practice. Recall from (\ref{hdefeq}), (\ref{Udefeq}) and (\ref{Vdefeq}), that
\[
h_{\ell}(t)=\frac{tV(t)}{U(t)}\ .
\]
By plugging the formal power series for $U(t)-1$ into $\frac{1}{1+x}=1-x+x^2-\cdots$, and then taking the product with the formal power series for $V(t)$, we immediately obtain the full asymptotic expansion
\[
h_{\ell}(t)=C_1 t+ C_2 t^2+C_3 t^3+\cdots
\]
where $C_1=1$, $C_2$, $C_3$, etc. are the constants which have been defined in \S\ref{constdefsec}.

Using explicit formulas for the reversion of formal power series, we have
the full asymptotic expansion for the inverse function
\[
h_{\ell}^{-1}(v)=D_1 v+ D_2 v^2+D_3 v^3+\cdots
\]
where $D_1=1$, $D_2$, $D_3$, etc. are the constants which have been defined in \S\ref{constdefsec}. We found it convenient to give a different name for the variable $v$ related to $t$ by $v=h_{\ell}(t)$.

\noindent{\bf Proof of the formula for the $D$ constants:}
Start from Theorem 1 and Corollary 1 from~\cite{AbdesselamSLC} which express the compositional inverse as a sum of tree-like Feynman diagrams. 
These trees have $p+1$ vertices of degree $1$: one root and $p$ leaves. They also have, for each $j\ge 2$,
$m_j$ (internal) vertices of degree $j+1$. The amplitude of such a tree graph is
\[
\mathcal{A}_{Fey}=v^{p}\times \prod_{j\ge 2}(-j! C_j)^{m_j}\ .
\]
Rather than deal directly with the complicated symmetry factor $\#Aut(E,\mathcal{F})$ in the denominator
of~\cite[Eq. (129)]{AbdesselamSLC}, it is better to use~\cite[Thm. 8]{AbdesselamSLC}
in order to forget about the bijection $\mathcal{C}$, i.e., use the forgetful map from Feynman diagram structures~\cite[Def. 12]{AbdesselamSLC} to pre-Feynman diagram structures~\cite[Def. 11]{AbdesselamSLC}. Now the symmetry factor in the denominator becomes the more manageable
$p!\prod_{j\ge 2}(m_j!\times j!^{m_j})$, see~\cite[Prop. 3]{AbdesselamSLC}.
However, one still has to sum over Wick contractions $\mathcal{C}$. For each Cayley tree, there are $\prod_{j\ge 2} j!^{m_j}$
such contractions, while the number of trees is
\[
\frac{(n-2)!}{(1-1)!^{1+p}\prod_{j\ge 2}\left((j+1)-1\right)!^{m_j}}
\]
where we made the vertex degrees apparent in the denominator.
Here, the total number of vertices is
\[
n=1+p+\sum_{j\ge 2}m_j\ .
\]
Also note that the quantity $\prod_{j\ge 2} j!^{m_j}$ appears twice in the numerator and twice in the denominator, and therefore cancels and is not visible in the final formulas.
Cleaning up the result of the above argument gives the explicit formula for the $D$ constants in \S\ref{constdefsec}.
\qed

We use the expansion of $h_{\ell}^{-1}$ to write
\[
(\ell-1)\varepsilon=
\frac{n}{k_n}h_{\ell}^{-1}\left((\ell-1)\frac{k_n}{n}\right)-\ell+1\sim
\sum_{j=2}^{\infty} D_j\times (\ell-1)^{j}\times \left(\frac{k_n}{n}\right)^{j-1}\ .
\]
as an asymptotic expansion when $n\rightarrow\infty$.
On the other hand, by inserting $\rho_n=\frac{n}{Z_{\ell}^{[\ell]}(t_n)}$, we have
\begin{eqnarray*}
-\ln(1+\eta_n)&=& -\ln\left[
\rho_n\times\frac{\mathcal{K}_{\ell}n^{\ell-1}}{(\ell-1)^{\ell}k_n^{\ell}} 
\right] \\
 & = & \ln\left[
\frac{(\ell-1)^{\ell}Z_{\ell}^{[\ell]}(t_n)}{\mathcal{K}_{\ell}}
\times\left(\frac{k_n}{n}\right)^{\ell}
\right] \\
 & = & \ln\left[
U(t_n)\times\left(\frac{h_{\ell}(t_n)}{t_n}\right)^{\ell}
\right]\ .
\end{eqnarray*}
Using
\[
\frac{h_{\ell}(t_n)}{t_n}\sim C_1+C_2 t_n+C_3 t_n^2+\cdots\ ,
\]
where $C_1=1$, and using 
\[
U(t_n)\sim A_0+A_1 t_n+A_2 t_n^2+\cdots
\]
where $A_0=1$, we obtain
\[
U(t_n)\times\left(\frac{h_{\ell}(t_n)}{t_n}\right)^{\ell}
\sim 1+E_1 t_n+E_2 t_n^2+\cdots
\]
by taking the product of the $\ell+1$ expansions involved.
The $E$ constants are as defined in \S\ref{constdefsec}.

We then insert into the last formula the asymptotic expansion
\[
t_n=h_{\ell}^{-1}\left((\ell-1)\frac{k_n}{n}\right)
\sim D_1 (\ell-1)\frac{k_n}{n}
+D_2 \left((\ell-1)\frac{k_n}{n}\right)^2+D_3\left((\ell-1)\frac{k_n}{n}\right)^3+\cdots
\]
which gives
\[
U(t_n)\times\left(\frac{h_{\ell}(t_n)}{t_n}\right)^{\ell}-1
\sim F_1\left(\frac{k_n}{n}\right)+F_2\left(\frac{k_n}{n}\right)^2+
F_3\left(\frac{k_n}{n}\right)^3+\cdots
\]
with the $F$ constants given in \S\ref{constdefsec}. We insert this in the expansion of $\ln(1+x)$, and deduce
\[
-\ln(1+\eta_n)\sim
G_1\left(\frac{k_n}{n}\right)+G_2\left(\frac{k_n}{n}\right)^2+
G_3\left(\frac{k_n}{n}\right)^3+\cdots
\]
with the $G$ constants given in \S\ref{constdefsec}.

Finally, adding the last expansion to that of $(\ell-1)\varepsilon_n$, we obtain
\[
\Delta_n\sim
H_1\left(\frac{k_n}{n}\right)+H_2\left(\frac{k_n}{n}\right)^2+
H_3\left(\frac{k_n}{n}\right)^3+\cdots
\]
with the $H$ constants defined in \S\ref{constdefsec}.
Hence, by keeping only the first $\ell-1$ terms explicit, we see that
\[
\Delta_n=\mathscr{E}_{\ell}\left(\frac{k_n}{n}\right)
+O\left(\left(\frac{k_n}{n}\right)^{\ell}\right)
\]
Since
\[
k_n\times\left(\frac{k_n}{n}\right)^{\ell}\sim s n^{-\frac{1}{\ell}}=o(1)\ ,
\]
Theorem \ref{mainthm1} is now established.
\qed

\section{The case of pairs and triples of permutations}

In the absence of an explicit reversion theorem for mixed asymptotic series with logarithms, we do the computation by hand to high enough order, essentially following the method of asymptotic iteration~\cite[Ch. 2]{deBruijn}.

\medskip
\noindent\underline{The $\ell=2$ case:}
Following similar steps as in the last section, we now have
\begin{eqnarray*}
\Delta_n & = & \varepsilon_n-\ln(1+\eta_n) \\
 & = & v_n^{-1}t_n-1+\ln U(t_n)+2\ln\left[t_n^{-1}h_2(t_n)\right]\ ,
\end{eqnarray*}
where we introduced the sequence $v_n:=\frac{k_n}{n}=h_2(t_n)$.
From (\ref{Z22eq}) and (\ref{Z21eq}), we derive
\begin{eqnarray*}
t^{-1}h_2 & = & \frac{1-\frac{\zeta(0)}{\zeta(2)}t\ln t+\frac{\zeta'(0)}{\zeta(2)}t+o(t)}{1+\frac{\zeta(0)}{\zeta(2)}t+o(t)} \\
 & = & \left(
1-\frac{\zeta(0)}{\zeta(2)}t\ln t+\frac{\zeta'(0)}{\zeta(2)}t+o(t)
\right)\times
\left(1+\frac{\zeta(0)}{\zeta(2)}t+o(t)\right) \\
 & = & 1-\frac{\zeta(0)}{\zeta(2)}t\ln t+\frac{\zeta'(0)-\zeta(0)}{\zeta(2)}t+o(t)\ .
\end{eqnarray*}
Hence,
\[
v_n=t_n-\frac{\zeta(0)}{\zeta(2)}t_n^2\ln t_n+\frac{\zeta'(0)-\zeta(0)}{\zeta(2)}t_n^2+o(t_n^2)\ ,
\]
which, turned around, gives
\[
t_n=v_n+\frac{\zeta(0)}{\zeta(2)}t_n^2\ln t_n+\frac{\zeta(0)-\zeta'(0)}{\zeta(2)}t_n^2+o(v_n^2)\ ,
\]
where we used $t_n\sim v_n$, to convert $o(t_n^2)$ into $o(v_n^2)$.
Let us introduce the quantity $\nu_n$ such that $t_n=v_n(1+\nu_n)$. We see from the above that $\nu_n=O(t_n\ln t_n)=O(v_n\ln v_n)$. By expanding the square, we easily see that
\[
t_n^2=v_n^2+o(v_n^2)\ .
\]
By expanding
\[
t_n^2\ln t_n=v_n^2(1+\nu_n)^2\left[\ln v_n+\ln(1+\nu_n)\right]\ ,
\]
we see that
\[
t_n^2\ln t_n=v_n^2\ln v_n+O(v_n^3(\ln v_n)^2)=v_n^2\ln v_n+o(v_n^2)\ .
\]
As a result,
\[
t_n=v_n+\frac{\zeta(0)}{\zeta(2)}v_n^2\ln v_n+\frac{\zeta(0)-\zeta'(0)}{\zeta(2)}v_n^2+o(v_n^2)\ ,
\]
and
\[
\varepsilon_n=v_n^{-1}t_n-1=
\frac{\zeta(0)}{\zeta(2)}v_n\ln v_n+\frac{\zeta(0)-\zeta'(0)}{\zeta(2)}v_n+o(v_n)\ .
\]
We have
\[
\ln U(t_n)=\ln\left[
1+\frac{\zeta(0)}{\zeta(2)}t_n+o(t_n)
\right]=\frac{\zeta(0)}{\zeta(2)}t_n+o(t_n)\ ,
\]
while
\begin{eqnarray*}
\ln\left[t_n^{-1} h_2(t_n)\right] & = & 
\ln\left[
1-\frac{\zeta(0)}{\zeta(2)}t_n\ln t_n+\frac{\zeta'(0)-\zeta(0)}{\zeta(2)}t_n+o(t_n)
\right] \\
 & = & -\frac{\zeta(0)}{\zeta(2)}t_n\ln t_n+\frac{\zeta'(0)-\zeta(0)}{\zeta(2)}t_n+o(t_n)\ .
\end{eqnarray*}
It is also easy to see, by expanding $v_n(1+\nu_n)\left[\ln v_n+\ln(1+\nu_n)\right]$, that
\[
t_n\ln t_n=v_n\ln v_n+O(v_n^2(\ln v_n)^2)=v_n\ln v_n+o(v_n)\ .
\]
Putting everything together we get
\[
\Delta_n=-\frac{\zeta(0)}{\zeta(2)}v_n\ln v_n+\frac{\zeta'(0)}{\zeta(2)}v_n+o(v_n)\ .
\]
The first two terms are what went into the definition of $\mathscr{E}_2$. The remainder $o(v_n)$, multiplied by $k_n$ of order  $n^{\frac{1}{2}}$ like $t_n^{-1}\sim v_n^{-1}$, then becomes $o(1)$ and Theorem \ref{mainthm2} is proved.
\qed

\medskip
\noindent\underline{The $\ell=3$ case:}
We proceed similarly to the $\ell=2$ case. We will spare the reader most of the tedious computations involved and will only present the highlights.
We again study
\begin{eqnarray*}
\Delta_n&=& 2\varepsilon_n-\ln(1+\eta_n) \\
 & = &  2\varepsilon_n+\ln U(t_n)+3\ln\left[t_n^{-1}h_3(t_n)\right]\ ,
\end{eqnarray*}
with
\[
\varepsilon=\frac{nt_n}{2k_n}-1\ .
\]
We perform the computations in terms of
$v_n=\frac{2k_n}{n}=h_3(t_n)$.

From (\ref{Zellasymeq}) with $\ell=3$, and from (\ref{Z32eq}),
we easily derive asymptotics for $h_3$ and thus $v_n$ in terms of $t_n$. Namely,
\begin{equation}
v_n=t_n+\frac{\zeta(0)}{2\zeta(3)}t_n^2-\frac{\zeta(-1)\zeta(0)}{\zeta(2)\zeta(3)}t_n^3\ln t_n
+C t_n^3+o(t_n^3)\ , 
\label{tnsolveeq}
\end{equation}
with the constant
\[
C:=\frac{1}{4\zeta(2)\zeta(3)^2}\left[
-\zeta(0)^2\zeta(2)
-2\zeta(-1)\zeta(0)\zeta(3)
+4\zeta'(-1)\zeta(0)\zeta(3)
+4\zeta(-1)\zeta'(0)\zeta(3)
\right]\ .
\]
We put all terms in (\ref{tnsolveeq}), other than the first $t_n$, 
on the other side and iterate this to eliminate $t_n$'s and trade them for $v_n$'s.
This time, a new phenomenon occurs, i.e., the conversion is not as straightforward as in the $\ell=2$ case.
For example, $t_n^2$ does not just give $v_n^2$, and we have a correction
\[
t_n^2=v_n^2-\frac{\zeta(0)}{\zeta(3)}v_n^3+o(v_n^3)\ .
\]
Likewise, in the computation of $-\ln(1+\eta_n)$, one needs
\begin{eqnarray*}
t_n & = & v_n-\frac{\zeta(0)}{2\zeta(3)}v_n^2+o(v_n^2) \\
t_n^2\ln t_n & = & v_n^2\ln v_n+o(v_n^2) \\
t_n^2 & = & v_n^2+o(v_n^2) \ .
\end{eqnarray*}
The rest is routine. We compute $\Delta_n$ up to $o(v_n^2)$ which multiplied by $k_n$ is $o(1)$. The explicit part goes into the definition of $\mathscr{E}_{3}$, and Theorem \ref{mainthm3} is established.
\qed

\section{Log-concavity}

We now provide a proof of Corollary \ref{maincor}.
Take a sequence $k_n$ as in the statements of Theorems \ref{mainthm1}, \ref{mainthm2}, or \ref{mainthm3}, whichever is appropriate given the value of $\ell$.
The shifted sequences $k_n-1$ and $k_n+1$ satisfy the same hypothesis of being asymptotically equivalent to $sn^{\frac{\ell-1}{\ell}}$, so we can apply the theorems to them too. Consider
\[
\Upsilon_n:=2\ln A(\ell,n,k_n)
-\ln A(\ell,n,k_n-1) -\ln A(\ell,n,k_n+1)\ .
\]
It can be written as
\[
\Upsilon_n=\ell\ \Upsilon_{n,{\rm main}}+\Upsilon_{n,{\rm err}}\ .
\]
Here $\Upsilon_{n,{\rm main}}$ comes from the only part of $\mathscr{S}_{\ell}$ which is not linear in $k$, i.e., $-k\ln k$, whereas  $\Upsilon_{n,{\rm err}}$ comes from $k\mathscr{E}_{\ell}\left(\frac{k}{n}\right)$.
We now compute
\begin{eqnarray*}
\Upsilon_{n,{\rm main}} & = & -2k_n\ln k_n
+(k_n-1)\ln(k_n-1)
+(k_n+1)\ln(k_n+1) \\
 & = & k_n\left[
\left(1-\frac{1}{k_n}\right)\ln\left(1-\frac{1}{k_n}\right)
+\left(1+\frac{1}{k_n}\right)\ln\left(1+\frac{1}{k_n}\right)
\right] \\
 & = & k_n\left[
\left(1-\frac{1}{k_n}\right)\left(-\frac{1}{k_n}-\frac{1}{2k_n^2}\right)+\left(1+\frac{1}{k_n}\right)\left(\frac{1}{k_n}-\frac{1}{2k_n^2}\right)
+o\left(\frac{1}{k_n^2}\right)
\right] \\
 & = & \frac{1}{k_n}+o\left(\frac{1}{k_n}\right) \\
 & = & s^{-1} n^{-\left(\frac{\ell-1}{\ell}\right)}+o\left( n^{-\left(\frac{\ell-1}{\ell}\right)}\right)\ .
\end{eqnarray*}
Now  the term $\Upsilon_{n,{\rm err}}$ is a linear combination of expressions of the form
\[
k_n\left(\frac{k_n}{n}\right)^a\left(\ln n-\ln k_n\right)^b
\]
with $a\ge 1$, which expanded becomes a combination of expressions
\[
n^{-a}(\ln n)^c k_n^{a+1}(\ln k_n)^d\ .
\]
Consider the function $f(k)=k^{a+1}(\ln k)^d$
and let us compute the following second order finite difference 
$L(k):=2f(k)-f(k-1)-f(k+1)$,
when $k$ is large. We claim that this finite difference is $O(k^{a-1}(\ln k)^d)$.
Indeed,
\begin{eqnarray*}
L(k) & = & 
2k^{a+1}(\ln k)^d
-(k-1)^{a+1}\left(\ln k+\ln\left(1-\frac{1}{k}\right)\right)^d \\
 & & -(k+1)^{a+1}\left(\ln k+\ln\left(1+\frac{1}{k}\right)\right)^d \\
 & = & k^{a+1}(\ln k)^d\left[2-\left(1-\frac{1}{k}\right)^{a+1}
-\left(1+\frac{1}{k}\right)^{a+1}
\right] \\
 & & -(k-1)^{a+1}\sum_{j=1}^{d}\binom{d}{j}(\ln k)^{d-j}\left(
\ln\left(1-\frac{1}{k}\right)\right)^j \\
 & & -(k+1)^{a+1}\sum_{j=1}^{d}\binom{d}{j}(\ln k)^{d-j}\left(
\ln\left(1+\frac{1}{k}\right)\right)^j\ .
\end{eqnarray*}
In the top line bracket the constant and linear term in $\frac{1}{k}$ cancel.
Hence that line's contribution is indeed $O(k^{a-1}(\ln k)^d)$.
The $j\ge 2$ terms in the sums also bring $O(k^{-2})$ factors multiplying $(k\pm 1)^{a+1}$. In the $j=1$ terms, one can expand $\ln(1+x)=x+O(x^2)$. The $O(x^2)$ part again brings $O(k^{-2})$ factors.
This leaves the linear part of $\ln(1+x)$ in the $j=1$ terms, i.e., the contribution
\[
\frac{1}{k}(\ln k)^{d-1}\times\sum_{i=0}^{a}(k-1)^{a-i}(k+1)^{i}
=O(k^{a-1}(\ln k)^d)\ ,
\]
and the claim is proved.
Now note that
\[
n^{-a}k_n^{a-1}\sim s^{a-1}\times n^{-\left(\frac{\ell-1}{\ell}\right)-\frac{a}{\ell}}=O(n^{-1})\ ,
\]
because $a\ge 1$.
Even after multiplication by powers of logarithms, this remains $o\left(n^{-\left(\frac{\ell-1}{\ell}\right)}\right)$. In conclusion, after putting everything together, we see that
\[
\Upsilon_n=\ell s^{-1} n^{-\left(\frac{\ell-1}{\ell}\right)}+o\left( n^{-\left(\frac{\ell-1}{\ell}\right)}\right)>0\ ,
\]
for $n$ large, 
and the corollary is proved.
\qed

\bigskip
\noindent{\bf Acknowledgements:}
{\small
We thank Shannon Starr for the collaboration~\cite{AbdesselamS}. We also thank him for seeing through the complexity of the numbers $A(\ell,n,k)$, and identifying their correct asymptotics, which led to the present work. We thank Ofir Gorodetsky, G\'erald Tenenbaum, and Mark Wildon for useful information and good advice on how to best handle the $Z$ function asymptotics.
} 


\begin{thebibliography}{999}

\bibitem{AbdesselamJPA}
A. Abdesselam, A physicist's proof of the Lagrange-Good multivariable inversion formula. J. Phys. A {\bf 36} (2003), no. 36, 9471--9477.

\bibitem{AbdesselamSLC}
A. Abdesselam, Feynman diagrams in algebraic combinatorics. S\'em. Lothar. Combin. {\bf 49} (2002/04), Art. B49c, 45 pp. (electronic).

\bibitem{AbdesselamCorr}
A. Abdesselam,
Non-Abelian correlation inequalities and stable determinantal polynomials.
Preprint arXiv:2207.07603[math-ph], 2022.

\bibitem{AbdesselamAC}
A. Abdesselam, Log-concavity with respect to the number of orbits for infinite tuples of commuting permutations. Ann. Comb. {\bf 29} (2025), no. 2, 563--573.

\bibitem{AbdesselamBDV}
A. Abdesselam, P. Brunialti, T. Doan, and P. Velie, 
A bijection for tuples of commuting permutations and a log-concavity conjecture.
Res. Number Theory {\bf 10} (2024), no. 2, Paper No. 45, 10 pp.

\bibitem{AbdesselamS}
A. Abdesselam, and S. Starr,
A central limit theorem for a generalization of the Ewens measure to random tuples of commuting permutations. Preprint arXiv:2505.11469[math.PR], 2025.

\bibitem{AndrewsAR}
G. E. Andrews, R. Askey, and R. Roy, {\it Special functions}. Encyclopedia of
Mathematics and its Applications {\bf 71}. Cambridge University Press, Cambridge, 1999.

\bibitem{Apostol}
T. M. Apostol,
{\it Introduction to Analytic Number Theory}.
Undergrad. Texts Math.,
Springer-Verlag, New York-Heidelberg, 1976.

\bibitem{BerndtE}
B. C. Berndt, and R. J. Evans,
Extensions of asymptotic expansions from Chapter 15 of Ramanujan's second notebook.
J. Reine Angew. Math. {\bf 361} (1985), 118--134.

\bibitem{BringmannFH}
K. Bringmann, J. Franke, and B. Heim,
Asymptotics of commuting $\ell$-tuples in symmetric groups and log-concavity.
Res. Number Theory {\bf 10} (2024), no. 4, Paper No. 83, 19 pp.

\bibitem{deBruijn}
N. G. de Bruijn,
{\it Asymptotic Methods in Analysis},
Second edition.
Bibl. Math., Vol. {\bf IV},
North-Holland Publishing Co., Amsterdam; P. Noordhoff Ltd., Groningen, 1961. 

\bibitem{BryanF}
J. Bryan, and J. Fulman,
Orbifold Euler characteristics and the number of commuting $m$-tuples in the symmetric groups.
Ann. Comb. {\bf 2} (1998), no. 1, 1--6.

\bibitem{DebruyneT}
G. Debruyne, and G. Tenenbaum,
The saddle-point method for general partition functions.
Indag. Math. (N.S.) {\bf 31} (2020), no. 4, 728--738.

\bibitem{Edwards}
H. M. Edwards,
{\it Riemann's Zeta Function}.
Reprint of the 1974 original (Academic Press, New York),
Dover Publications, Inc., Mineola, NY, 2001. 

\bibitem{FlajoletGD}
P. Flajolet, X. Gourdon, and P. Dumas,
Mellin transforms and asymptotics: harmonic sums.
Special volume on mathematical analysis of algorithms
Theoret. Comput. Sci.  {\bf 144} (1995), no. 1--2, 3--58.

\bibitem{HeimN}
B. Heim, and M. Neuhauser,
Horizontal and vertical log-concavity.
Res. Number Theory {\bf 7} (2021), no. 1, Paper No. 18, 12 pp.

\bibitem{Hwang}
H.-K. Hwang,
Limit theorems for the number of summands in integer partitions.
J. Combin. Theory Ser. A {\bf 96} (2001), no. 1, 89--126.

\bibitem{Kowalski}
E. Kowalski,
{\it Un Cours de Th\'eorie Analytique des Nombres}.
Cours Sp\'ec., {\bf 13},
Soci\'et\'e Math\'ematique de France, Paris, 2004. 

\bibitem{Lagarias}
J. C. Lagarias,
Euler's constant: Euler's work and modern developments.
Bull. American Math. Soc. (N.S.) {\bf 50} (2013), no. 4, 527--628.

\bibitem{LipnikMT1}
G. F. Lipnik, M. G. Madritsch, and R. F. Tichy,
A central limit theorem for integer partitions into small powers.
Monatsh. Math. {\bf 203} (2024), no. 1, 149--173.

\bibitem{LipnikMT2}
G. F. Lipnik, M. G. Madritsch, and R. F. Tichy,
A local limit theorem for integer partitions into small powers.
Preprint arXiv:2311.09203 [math.CO], 2023.

\bibitem{MadritschW}
M. Madritsch, and S. Wagner,
A central limit theorem for integer partitions.
Monatsh. Math. {\bf 161} (2010), no. 1, 85--114.

\bibitem{MoserW}
L. Moser, M. Wyman,
Asymptotic development of the Stirling numbers of the first kind.
J. London Math. Soc. {\bf 33} (1958), 133--146.

\bibitem{NinhamHFG}
B. W. Ninham, B. D. Hughes, N. E. Frankel, and M. L. Glasser,
M\"{o}bius, Mellin, and mathematical physics.
Physica A {\bf 186} (1992), no. 3-4, 441--481.

\bibitem{PemantleW}
R. Pemantle, and M. C. Wilson,
{\it Analytic Combinatorics in Several Variables}.
Cambridge Stud. Adv. Math. {\bf 140},
Cambridge University Press, Cambridge, 2013.

\bibitem{Riemann}
B. Riemann, \"{U}ber die Anzahl der Primzahlen unter einer gegebenen Gr\"{o}sse. Monatsberichte
der Berliner Akademie (1859) 671--680.

\bibitem{Starr}
S. Starr, Some observations about the ``generalized abundancy index''.
Preprint arXiv:2505.07051[math.CO], 2025.

\bibitem{Tao}
T. Tao,
{\it Compactness and Contradiction}.
American Mathematical Society, Providence, RI, 2013. 

\bibitem{TenenbaumBook}
G. Tenenbaum,
{\it Introduction to Analytic and Probabilistic Number Theory}.
Third edition. Translated from the 2008 French edition by Patrick D. F. Ion.
Grad. Stud. Math. {\bf 163}
American Math. Soc., Providence, RI, 2015.

\bibitem{TenenbaumWL}
G. Tenenbaum, J. Wu, and Y.-L. Li,
Power partitions and saddle-point method.
J. Number Theory {\bf 204} (2019), 435--445.

\bibitem{Titchmarsh}
E. C. Titchmarsh,
{\it The Theory of The Riemann Zeta-Function}.
Second edition. Edited and with a preface by D. R. Heath-Brown.
The Clarendon Press, Oxford University Press, New York, 1986. 

\bibitem{Wigert}
S. Wigert,
Sur la s\'erie de Lambert et son application \`a la th\'eorie des nombres.
Acta Math. {\bf 41} (1916), no. 1, 197--218.

\bibitem{Zagier}
D. Zagier, The Mellin transfom and related analytic techniques. 
Appendix to Ch. 6 of the book by E. Zeidler {\it Quantum Field Theory I: Basics in Mathematics and Physics. A Bridge Between Mathematicians and Physicists}, pp. 305-323.
Springer-Verlag, Berlin-Heidelberg-New York, 2006. 

\end{thebibliography}
\end{document}